\documentclass[leqno,10pt]{amsart}
\usepackage{amsmath,amscd}
\usepackage{epsfig,graphics,texdraw}
\usepackage{amssymb,latexsym}
\usepackage{mathrsfs}
\usepackage{eucal}
\newtheorem{thm}{\bf Theorem}[section]
\newtheorem{df}[thm]{\bf Definition}
\newtheorem{prop}[thm]{\bf Proposition}
\newtheorem{cor}[thm]{\bf Corollary}
\newtheorem{lem}[thm]{\bf Lemma}
\newtheorem{rem}[thm]{\bf Remark}
\newtheorem{ex}[thm]{\bf Example}
\newcommand{\bs}{\boldsymbol}

\newcommand{\B}{\mathbf{B}}
\newcommand{\cS}{\mathcal{S}}
\newcommand{\C}{\mathscr{C}}
\newcommand{\W}{\mathbf{W}}
\newcommand{\cP}{\mathscr{P}}

\newcommand{\pf}{\noindent{\bfseries Proof. }}

\newcommand{\bi}{\bs{\rm i}}
\newcommand{\bj}{\bs{\rm j}}

\newcommand{\bk}{\bs{\rm k}}
\newcommand{\bl}{\bs{\rm l}}

\newcommand{\M}{{\bf M}}

\newcommand{\hf}{\frac{1}{2}}
\newcommand{\N}{\mathcal{N}}

\numberwithin{equation}{section}

\begin{document}
\title[Crystals and quasi-symmetric functions]
{Crystal graphs for general linear Lie superalgebras and
quasi-symmetric functions}
\author{JAE-HOON KWON}
\address{Department of Mathematics \\ University of Seoul \\ 90
Cheonnong-dong, Dongdaemun-gu \\ Seoul 130-743, Korea }
\email{jhkwon@uos.ac.kr }

\thanks{This research was supported by KRF Grant $\sharp$2005-070-C00004} 

\begin{abstract}
We give a new representation theoretic interpretation of the ring of
quasi-symmetric functions. This is obtained by showing that the
super analogue of the Gessel's fundamental quasi-symmetric function
can be realized as the character of an irreducible crystal for the
Lie superalgebra $\frak{gl}_{n|n}$ associated to its non-standard
Borel subalgebra with a maximal number of odd isotropic simple
roots. We also present an algebraic characterization of these super
quasi-symmetric functions.
\end{abstract}

\maketitle


\section{Introduction}

The notion of a quasi-symmetric function is a generalization of a
symmetric function introduced by Gessel \cite{G}. The ring of
quasi-symmetric functions $QSym$ has many interesting features, for
example, it is a Hopf algebra whose dual is isomorphic to the ring
of non-commutative symmetric functions introduced by Gefand et al.
\cite{GDLLRT}. As in the case of symmetric functions, $QSym$ and its
dual have nice representation theoretic interpretations related to
the representations of degenerate Hecke algebras and quantum groups
of type ${\rm A}$ \cite{KT97,KT99}. In \cite{KT99}, Krob and Thibon
showed that Gessel's fundamental quasi-symmetric functions in $n$
variables, which form a basis of $QSym$ in $n$ variables, are the
characters of the irreducible polynomial representations of the
degenerate quantum group associated to $\frak{gl}_n$. Moreover, the
action of the degenerate quantum group determines a quasi-crystal
graph structure on the irreducible polynomial representation, which
is a subgraph of the Kashiwara's crystal graph of Young tableaux
(see also \cite{Th01} for a general review on this topic).

The purpose of the present paper is to give a new representation
theoretic interpretation of $QSym$ using the crystal base theory for
contragredient Lie superalgebras developed by Benkart, Kang, and
Kashiwara \cite{BKK}. This work was motivated by observing that
unlike finite dimensional complex simple Lie algebras and their
crystals, two Borel subalgebras of a general linear Lie superalgebra
$\frak{gl}_{m|n}$ are not necessarily conjugate to each other, and
the corresponding $\frak{gl}_{m|n}$-crystal structures are different
in general since they may have different systems of simple roots.
Among the Borel conjugacy classes of $\frak{gl}_{n|n}$, we study in
detail a crystal structure when the simple roots associated to a
Borel subalgebra are all odd and hence isotropic, while the crystals
associated to a standard Borel subalgebra of $\frak{gl}_{m|n}$ with
a single odd simple root is discussed explicitly in \cite{BKK}. Then
it turns out that the combinatorics of these non-standard crystals
is remarkably simple and nice. Our main result is that as a
$\frak{gl}_{n|n}$-crystal associated to this non-standard Borel
conjugacy class, the connected components in the crystal of tensor
powers of the natural representation are parameterized by
compositions, and each of them can be realized as the set of (super)
quasi-ribbon tableaux of a given composition shape
(cf.~\cite{KT99}). An insertion algorithm for these non-standard
crystals can be derived in a standard way, and it is shown to be
compatible with our crystal structure. Hence, we obtain an explicit
description of the Knuth equivalence or crystal equivalence for this
case. Furthermore, we show that the set of Stembridge's enriched
$P$-partitions \cite{St97} is naturally equipped with this
non-standard crystal structure, and then apply its combinatorial
properties to decomposition of crystals including tensor product
decompositions. As a corollary, it follows immediately that the
irreducible characters of these non-standard
$\frak{gl}_{n|n}$-crystals are equal to the super analogue of
fundamental quasi-symmetric functions, which are defined by a
standard method of superization \cite{HHRU}, and they form a ring
isomorphic to $QSym$ under a suitable limit.

One may regard the notion of non-standard $\frak{gl}_{n|n}$-crystals
as the counterpart of the Krob and Thibon's quasi-crystals, and
hence as a quasi-analogue of standard $\frak{gl}_{n|n}$-crystals. In
this sense, one can define a quasi-analogue of a standard
$\frak{gl}_{m|n}$-crystal for arbitrary $m$ and $n$ by a crystal
associated to its Borel conjugacy class having a maximal number of
odd isotropic simple roots. This enables us to explain the relation
between symmetric functions and quasi-symmetric functions as a
special case of branching rules between $\frak{gl}_{m|n}$-crystals
associated to a standard and a non-standard Borel conjugacy classes.
We discuss a crystal structure associated to a Borel conjugacy class
having a maximal number of odd isotropic simple roots, and classify
its connected components occurring in the crystal of tensor powers
of the natural representation, which are parameterized by certain
pairs of a partition and a composition. Then we obtain an explicit
branching decomposition of standard $\frak{gl}_{m|n}$-crystals into
non-standard ones.

The paper is organized as follows. In Section 2, we introduce the
notion of abstract crystals for general linear Lie superalgebras and
recall some basic properties. In Section 3, we discuss in detail a
crystal structure associated to a Borel subalgebra whose simple
roots are all odd isotropic. In Section 4, we study a non-standard
crystal structure on the set of enriched $P$-partitions and its
applications. In Section 5, we consider in general a crystal
structure associated to a Borel conjugacy class with a maximal
number of odd isotropic simple roots. Finally, in Section 6, we give
an algebraic characterization of irreducible characters of
non-standard crystals discussed in the previous section.\vskip 5mm

\section{Crystal graphs for general linear Lie superalgebras}

\subsection{Lie superalgebra $\frak{gl}_{\cS}$ and crystal graphs}
Suppose that ${\cS}$ is a $\mathbb{Z}_2$-graded set with a linear
ordering $\prec$. Let $\mathbb{C}^{\mathcal S}$ be the complex
superspace with a basis $\{\,v_b\,|\,b\in \cS\,\}$. Let
$\frak{gl}_{\cS}=\frak{gl}(\mathbb{C}^{\cS})$ denote the Lie
superalgebra of complex linear transformations on $\mathbb{C}^{S}$
that vanish on a subspace of finite codimension. We call
$\frak{gl}_{{\cS}}$ a {\em general linear Lie superalgebra}
(cf.~\cite{Kac}). In particular, when
${\cS}=\{\,k\in\mathbb{Z}^{\times}\,|\,-m\leq k\leq n\,\}=[m|n]$
($m,n\geq 0$) with a usual linear ordering,
${\cS}_0=\{-m,\ldots,-1\}$ and ${\cS}_1=\{1,\ldots,n\}$,
$\frak{gl}_{\cS}$ is often denoted by $\frak{gl}_{m|n}$. We may
identify $\frak{gl}_{\cS}$ with the space of complex matrices
generated by the elementary matrices $e_{ab}$ ($a,b\in {\cS}$). Then
$\frak{h}=\sum_{b\in {\cS}}\mathbb{C}e_{bb}$ is a Cartan subalgebra.
Let $\langle\ , \ \rangle$ be the natural pairing on
$\frak{h}^*\times \frak{h}$. Let $\epsilon_a\in\frak{h}^*$ be
determined by $\langle \epsilon_a,e_{bb} \rangle=\delta_{ab}$ for
$a,b\in {\cS}$. Let $P=\bigoplus_{b\in {\cS}} \mathbb{Z}
\epsilon_{b}$ be the weight lattice for $\frak{gl}_{\cS}$, which is
a free abelian group generated by $\epsilon_b$ ($b\in {\cS}$). There
is a natural symmetric $\mathbb{Z}$-bilinear form $(\ ,\,)$ on $P$
given by $(\epsilon_a,\epsilon_{b})=(-1)^{|a|}\delta_{ab}$ for
$a,b\in {\cS}$, where $|a|$ denotes the degree of $a$.

Let $\frak{b}$ be the subalgebra of $\frak{gl}_{\cS}$ spanned by
$e_{ab}$ for $a\preccurlyeq b\in{\cS}$, which we call a Borel
subalgebra. The set of positive simple roots and the set of positive
roots of $\frak{b}$ with respect to $\frak{h}$ are given by
\begin{equation}\label{simple roots}
\begin{split}
\Delta\ (\text{or }\Delta_{\cS})&=\{\,\epsilon_{b}-\epsilon_{b'}\
\text{for
a successive pair $b\prec b'\in {\cS}$}\,\}, \\
\Phi\ (\text{or }\Phi_{\cS})&=\{\epsilon_{b}-\epsilon_{b'}\
\text{for
 $b\prec b'\in {\cS}$}\,\}.
\end{split}
\end{equation}
Note that $(\alpha,\alpha)=\pm 2, 0$ for $\alpha \in \Phi$. We put
$\ell_{\alpha}=(-1)^{|b|}$ for
$\alpha=\epsilon_{b}-\epsilon_{{b'}}\in\Delta$. Let
$Q=\bigoplus_{\alpha \in \Delta} \mathbb{Z} \alpha$ be the root
lattice of $\frak{gl}_{{\cS}}$. A partial ordering on $P$ with
respect to $Q$ is given by $\lambda\geq \mu$ if and only if
$\lambda-\mu\in \sum_{\alpha \in \Delta} \mathbb{Z}_{\geq 0} \alpha$
for $\lambda,\mu\in P$.

Now, let us introduce the notion of abstract crystal graphs for
$\frak{gl}_{{\cS}}$. Our definition is based on the crystal base
theory for the quantized enveloping algebra of a contragredient Lie
superalgebra  developed by Benkart, Kang, and Kashiwara \cite{BKK}.
Throughout the paper, ${\bf 0}$ denotes a formal symbol.

\begin{df}\label{crystal graph}{\rm   A {\it crystal graph for
$\frak{gl}_{{\cS}}$}, or simply {\it $\frak{gl}_{{\cS}}$-crystal},
is a set $B$ together with the maps ${\rm wt}  : B \rightarrow P$,
$\varepsilon_{\alpha}, \varphi_{\alpha}: B \rightarrow
\mathbb{Z}_{\geq 0}$, ${e}_{\alpha}, {f}_{\alpha}: B \rightarrow
B\cup\{{\bf 0}\}$ (${\alpha}\in \Delta$) such that for $b,b'\in B$,
\begin{itemize}

\item[(1)] if $\alpha$ is isotropic (that is,  $(\alpha,\alpha)=0$), then
$$\varphi_{\alpha}(b)+ \varepsilon_{\alpha}(b)=
\begin{cases}
1, & \text{if $({\rm wt}(b),\alpha)\neq 0$}, \\
0, & \text{otherwise},
\end{cases}$$

\item[(2)] if $\alpha$ is non-isotropic, then
$\varphi_{\alpha}(b)-\varepsilon_{\alpha}(b)=\ell_{\alpha}({\rm
wt}(b),\alpha)$,

\item[(3)] if ${e}_{\alpha} b \in B$, then
$\varepsilon_{\alpha}({e}_{\alpha} b) = \varepsilon_{\alpha}(b) -
1$, $\varphi_{\alpha}({e}_{\alpha} b) = \varphi_{\alpha}(b) + 1$,
and ${\rm wt}(e_{\alpha}b)={\rm wt}(b)+\alpha$,

\item[(4)] if ${f}_{\alpha} b \in B$, then
$\varepsilon_{\alpha}({f}_{\alpha} b) = \varepsilon_{\alpha}(b) +
1$, $\varphi_{\alpha}({f}_{\alpha} b) = \varphi_{\alpha}(b) - 1$,
and ${\rm wt}(f_{\alpha}b)={\rm wt}(b)-\alpha$,

\item[(5)] ${f}_{\alpha} b = b'$ if and only if $b = {e}_{\alpha}
b'$.
\end{itemize}}
\end{df}\vskip 3mm
A $\frak{gl}_{{\cS}}$-crystal $B$ becomes a  $\Delta$-colored
oriented graph, where $b\stackrel{\alpha}{\rightarrow}b'$ if and
only if $b'=f_{\alpha}b$ $(\alpha\in \Delta)$.  We call
$e_{\alpha}$, $f_{\alpha}$ the Kashiwara operators.

\begin{rem}{\rm
(1) If $(\alpha,\alpha)=2$ for all $\alpha\in \Delta$, or
${\cS}_0={\cS}$, then $\Phi$ is equal to the set of positive roots
of type ${\rm A}_{n-1}$ with $n=|{\cS}|$, and
$\frak{gl}_{{\cS}}$-crystals are the crystal graphs for ${\rm
A}_{n-1}$ (cf.~\cite{Kas90,Kas94}). On the other hand, if
$(\alpha,\alpha)=-2$ (equivalently, $\alpha$ is non-isotropic and
$\ell_{\alpha}=-1$) for all $\alpha\in \Delta$, then $\Phi$ can be
identified with the set of negative roots of type ${\rm A}_{n-1}$,
and $\frak{gl}_{{\cS}}$-crystals are dual crystal graphs for ${\rm
A}_{n-1}$.

(2) For an isotropic simple root $\alpha$, we have
$\varepsilon_{\alpha}(b)=1$ (resp. $\varphi_{\alpha}(b)=1$)  if
$e_{\alpha}b\neq {\bf 0}$ (resp. $f_{\alpha}b\neq {\bf 0}$), and
$e_{\alpha}^2 b=f_{\alpha}^2 b={\bf 0}$ (cf.~\cite{BKK}). }
\end{rem}

Suppose that $B_1$ and $B_2$ are $\frak{gl}_{{\cS}}$-crystals. We
say that $B_1$ is isomorphic to $B_2$, and write $B_1\simeq B_2$ if
there is an isomorphism of $\Delta$-colored oriented graphs which
preserves ${\rm wt}$, $\varepsilon_{\alpha}$, and $\varphi_{\alpha}$
($\alpha\in \Delta$). For $b_i\in B_i$ ($i=1,2$), let $C(b_i)$
denote the connected component of $b_i$ as an $\Delta$-colored
oriented graph. We say that $b_1$ is {\it
$\frak{gl}_{{\cS}}$-equivalent to} $b_2$, if there is an isomorphism
of crystal graphs $C(b_1)\rightarrow C(b_2)$ sending $b_1$ to $b_2$,
and write $b_1 {\simeq} b_2$ for short.

We define the tensor product $B_1\otimes B_2=\{\,b_1\otimes
b_2\,|\,b_i\in B_i\,\, (i=1,2)\,\}$ as follows;

{\allowdisplaybreaks
\begin{equation*}
{\rm wt}(b_1\otimes b_2)={\rm wt}(b_1)+{\rm wt}(b_2),
\end{equation*}

(1) for an isotropic simple root $\alpha\in \Delta$,
\begin{equation*}
\begin{split}
&\chi_{\alpha}(b_1\otimes b_2)= \\
&\begin{cases} \chi_{\alpha}(b_1), & \text{if ($\ell_{\alpha}=1$,
$(\alpha,{\rm wt}(b_1))\neq 0$)
or ($\ell_{\alpha}=-1$, $(\alpha,{\rm wt}(b_2))=0$)}, \\
\chi_{\alpha}(b_2), & \text{if ($\ell_{\alpha}=-1$, $(\alpha,{\rm
wt}(b_2))\neq 0$) or ($\ell_{\alpha}=1$, $(\alpha,{\rm
wt}(b_1))=0$)},
\end{cases}
\end{split}
\end{equation*}
\begin{equation*}
\begin{split}
&x_{\alpha}(b_1\otimes b_2)= \\
&\begin{cases} {x}_{\alpha} b_1 \otimes b_2, & \text{if
($\ell_{\alpha}=1$, $(\alpha,{\rm wt}(b_1))\neq 0$)
or ($\ell_{\alpha}=-1$, $(\alpha,{\rm wt}(b_2))=0$)}, \\
b_1 \otimes {x}_{\alpha} b_2, & \text{if ($\ell_{\alpha}=-1$,
$(\alpha,{\rm wt}(b_2))\neq 0$) or ($\ell_{\alpha}=1$, $(\alpha,{\rm
wt}(b_1))=0$)},
\end{cases}
\end{split}
\end{equation*}
where $\chi=\varepsilon, \varphi$ and $x=e,f$,

(2) for a non-isotropic simple root $\alpha\in\Delta$,
\begin{equation*}
\begin{split}
&\chi_{\alpha}(b_1\otimes b_2)= \\
&\begin{cases} {\rm
max}\{\chi_{\alpha}(b_1),\chi_{\alpha}(b_2)-\ell_{\alpha}(\alpha,{\rm
wt}(b_1))\}, & \text{if ($\chi=\varepsilon$, $\ell_{\alpha}=1$)
or ($\chi=\varphi$, $\ell_{\alpha}=-1$)}, \\
{\rm max}\{\chi_{\alpha}(b_1)+\ell_{\alpha}(\alpha,{\rm
wt}(b_2)),\chi_{\alpha}(b_2)\}, & \text{if ($\chi=\varphi$,
$\ell_{\alpha}=1$) or ($\chi=\varepsilon$, $\ell_{\alpha}=-1$)},
\end{cases}
\end{split}
\end{equation*}
\begin{equation*}
\begin{split}
&e_{\alpha}(b_1\otimes b_2)= \\
&\begin{cases} {e}_{\alpha} b_1 \otimes b_2, & \text{if
($\ell_{\alpha}=1$, $\varphi_{\alpha}(b_1)\geq
\varepsilon_{\alpha}(b_2)$)
or ($\ell_{\alpha}=-1$, $\varphi_{\alpha}(b_2)< \varepsilon_{\alpha}(b_1)$)}, \\
b_1 \otimes {e}_{\alpha} b_2, & \text{if ($\ell_{\alpha}=1$,
$\varphi_{\alpha}(b_1)< \varepsilon_{\alpha}(b_2)$) or
($\ell_{\alpha}=-1$, $\varphi_{\alpha}(b_2)\geq
\varepsilon_{\alpha}(b_1)$)},
\end{cases}
\end{split}
\end{equation*}
\begin{equation*}
\begin{split}
&f_{\alpha}(b_1\otimes b_2)= \\
&\begin{cases} {f}_{\alpha} b_1 \otimes b_2, & \text{if
($\ell_{\alpha}=1$, $\varphi_{\alpha}(b_1)>
\varepsilon_{\alpha}(b_2)$)
or ($\ell_{\alpha}=-1$, $\varphi_{\alpha}(b_2)\leq \varepsilon_{\alpha}(b_1)$)}, \\
b_1 \otimes {f}_{\alpha} b_2, & \text{if ($\ell_{\alpha}=1$,
$\varphi_{\alpha}(b_1)\leq \varepsilon_{\alpha}(b_2)$) or
($\ell_{\alpha}=-1$, $\varphi_{\alpha}(b_2)>
\varepsilon_{\alpha}(b_1)$)},
\end{cases}
\end{split}
\end{equation*}
where $\chi=\varepsilon, \varphi$. We assume that ${\bf 0}\otimes
b_2=b_1\otimes {\bf 0}={\bf 0}$.} It is straightforward to check
that $B_1\otimes B_2$ is a $\frak{gl}_{{\cS}}$-crystal.

\subsection{Crystals of semistandard tableaux} We may regard ${\cS}$ as a
$\frak{gl}_{{\cS}}$-crystal associated to the natural representation
$\mathbb{C}^{\cS}$, where $b\stackrel{\alpha}{\rightarrow}b'$ for a
successive pair $b\prec b'$ in ${\cS}$ with
$\alpha=\epsilon_b-\epsilon_{b'}$, ${\rm wt}(b)=\epsilon_b$, and
$\varepsilon_{\beta}(b)$ (resp. $\varphi_{\beta}(b)$) is the number
of $\beta$-colored arrows coming into $b$ (resp. going out of $b$)
for $\beta\in\Delta$.

Let $\W=\W_{{\cS}}^{}$ be the set of words of finite length with
alphabets in ${\cS}^{}$. The empty word is denoted by $\emptyset$.
Then $\W$ is a $\frak{gl}_{{\cS}}^{}$-crystal since we may view each
non-empty word $w=w_1\cdots w_r$ as $w_1\otimes\cdots\otimes w_r\in
{\cS}^{\otimes r}$, where $\{\emptyset\}$ forms a trivial crystal
graph, that is, ${\rm wt}(\emptyset)=0$,
$e_{\alpha}\emptyset=f_{\alpha}\emptyset={\bf 0}$, and
$\varepsilon_{\alpha}(\emptyset)=\varphi_{\alpha}(\emptyset)=0$ for
$\alpha\in\Delta$.

Let $\cP$ be the set of partitions. As usual, we identify
$\lambda=(\lambda_k)_{k\geq 1}\in\cP$ with its Young diagram. For
non-negative integers $m$ and $n$, we denote by $\cP_{m|n}$ the set
of all $(m,n)$-hook partitions $\lambda$, that is,
$\lambda_{m+1}\leq n$.

For $\lambda\in\cP$, let $\B(\lambda)=\B_{{\cS}}(\lambda)$ be the
set of {\it semistandard tableaux of shape $\lambda$}, that is,
tableaux $T$ obtained by filling a Young diagram $\lambda$ with
entries in ${\cS}^{}$ such that (1) the entries in each row (resp.
column) are weakly increasing from left to right (resp. from top to
bottom), and (2) the entries in ${\cS}^{}_0$ (resp. ${\cS}^{}_1$)
are strictly increasing in each column (resp. row) (cf.~\cite{BR}).
Note that $\B(\lambda)$ is non-empty if and only if
$\lambda\in\cP_{m|n}$, where $|{\cS}_0|=m$ and $|{\cS}_1|=n$. We
assume that $\cP_{m|n}=\cP$ if either $m$ or $n$ is infinite.

For $T\in\B(\lambda)$, $T(i,j)$ denotes the entry of $T$ located in
the $i$th row from the top and the $j$th column from the left. Then
an {\it admissible reading} is an embedding $\psi : \B(\lambda)
\rightarrow \W$ given by reading the entries of $T$ in $\B(\lambda)$
in such a way that $T(i,j)$ should be read before $T(i+1,j)$ and
$T(i,j-1)$ for each $i,j$. Then, by similar arguments as in
\cite{BKK}, we can check that the image of $\B(\lambda)$ under
$\psi$ together with ${\bf 0}$ is stable under
$e_{\alpha},f_{\alpha}$ ($\alpha\in \Delta$), and hence
$\B(\lambda)$ is a subcrystal of $\W$ (cf.~\cite{Kas94}), which does
not depend on the choice of $\psi$.

\begin{prop}[cf.~\cite{BKK}]\label{BKK} For $\lambda\in\cP$,
$\B(\lambda)$ is a $\frak{gl}_{{\cS}}$-crystal.
\end{prop}

For $b\in {\cS}$ and $T\in\B(\lambda)$, let us denote by $b
\Rightarrow T$ the tableau in $\B(\mu)$ obtained by { Schensted's
column bumping insertion} (cf.~\cite{BR,Rem}), where $\mu\in\cP$ is
given by adding a node at $\lambda$. Now, for a given word
$w=w_1\cdots w_r\in\W$, we define its $P$-tableau by
\begin{equation}\label{column insertion}
{\bf P}(w)=w_r \Rightarrow(\cdots\Rightarrow(w_2 \Rightarrow w_1)).
\end{equation}

\begin{prop}[cf.~\cite{BKK}]
For $w\in\W$, we have $w\simeq {\bf P}(w)$. 
\end{prop}

Let us consider the case when ${\cS}=[m|n]$. For
$\lambda\in\cP_{m|n}$, let $H_{\lambda}$ be the unique element in
$\B_{[m|n]}(\lambda)$ such that ${\rm
wt}(H_{\lambda})=\sum_{i=1}^m\lambda_i\epsilon_{-m+i-1}+\sum_{j=
1}^{n}\mu_j\epsilon_j$, where $\mu=(\mu_j)$ is the transpose of
$(\lambda_{m+1},\lambda_{m+2},\ldots)$. Then the following is one of
the main results in \cite{BKK}.

\begin{thm}[\cite{BKK}]\label{Bmnlambda}
For $\lambda\in \cP_{m|n}$, $\B_{[m|n]}(\lambda)$ is connected with
a unique highest weight element $H_{\lambda}$. 
\end{thm}

\begin{rem}{\rm
(1) One of the most important properties of the crystal graph
$\B_{[m|n]}(\lambda)$ is the existence of fake highest weight
vectors. That is, $\B_{[m|n]}(\lambda)$ can have an element $T$ such
that $e_{\alpha}T={\bf 0}$ for all $\alpha\in\Delta_{[m|n]}$, but
$T\neq H_{\lambda}$.

(2) For a general ${\cS}$, $\B_{{\cS}}(\lambda)$ is not necessarily
connected. We will see some examples in later sections.}
\end{rem}

\section{Crystals of quasi-ribbon tableaux}
In this section, we discuss in detail the crystal graphs for
$\frak{gl}_{\N}$, where
$$\N=\text{\tiny $\hf$}\mathbb{Z}_{>0}=\Big{\{} \,\text{\tiny$\frac{1}{2}$} \prec 1 \prec \text{\tiny$\frac{3}{2}$} \prec 2 \prec \text{\tiny$\frac{5}{2}$} \prec \cdots\,\Big{\}},$$
with $\N_0=\mathbb{Z}_{>0}$ and $\N_1=\hf+\mathbb{Z}_{\geq 0}$. Then
$\Delta_{\N}=\{\,\alpha_{r}=\epsilon_{r}-\epsilon_{r+\hf},\
(r\in\hf\mathbb{Z}_{>0}) \,\}$. Denote by $I=\hf\mathbb{Z}_{>0}$ the
index set for the simple roots. The associated Dynkin diagram is
\vskip 5mm
\begin{center}
\includegraphics{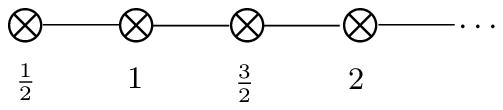}\ \ \ .
\end{center}
Note that $\ell_{\alpha_i}=(-1)^{2i}$  for $i\in I$.

\subsection{Quasi-ribbon tableaux}
The $\frak{gl}_{\N}$-crystal $\N$ is given by
\begin{equation*}
\text{\tiny$\hf$}\ \stackrel{\hf}{\longrightarrow}\ 1\
\stackrel{1}{\longrightarrow} \text{\tiny$\frac{3}{2}$}
\stackrel{\frac{3}{2}}{\longrightarrow}\ 2\
\stackrel{2}{\longrightarrow}\ \text{\tiny$\frac{5}{2}$} \
{\longrightarrow} \cdots\ \ .
\end{equation*}
From the tensor product rule of crystals given in Section 2.1, the
actions of $e_i=e_{\alpha_i}, f_i=f_{\alpha_i}$ on $w=w_1\cdots
w_r\in \W=\W_{\N}$ ($i\in I$) can be described more explicitly;
\begin{itemize}
\item[(1)] for $i\in\mathbb{Z}_{>0}$, choose the smallest $k$ ($1\leq k\leq
r$) such that $(\alpha_i,{\rm wt}(w_k))\neq 0$. Then $e_i w$ (resp.
$f_i w$) is obtained by applying $e_i$ (resp. $f_i$) to $w_k$. If
there is no such $k$, then   $e_i w= {\bf 0}$ (resp. $f_iw={\bf
0}$).

\item[(2)] for $i\in\hf+\mathbb{Z}_{\geq 0}$, choose the largest $k$ ($1\leq k\leq
r$) such that $(\alpha_i,{\rm wt}(w_k))\neq 0$. Then $e_i w$ (resp.
$f_i w$) is obtained by applying $e_i$ (resp. $f_i$) to $w_k$. If
there is no such $k$, then   $e_i w= {\bf 0}$ (resp. $f_iw={\bf
0}$).
\end{itemize}
Note that $\varepsilon_i(w)$ (resp. $\varphi_i(w)$) is the number of
$i$-colored arrows coming into $w$ (resp. going out of $w$) for
$w\in \W$ and $i\in I$.

A {\it composition} is a finite sequence of positive integers
$\alpha=(\alpha_1,\ldots,\alpha_t)$ for some $t\geq 1$, and we write
$\alpha\vDash r$ if $\sum_{i=1}^t\alpha_i=r$. Let
$S(\alpha)=\{\,\alpha_1+\cdots+\alpha_i\,|\,i=1,\ldots,t-1\,\}$,
which is a subset of $\{1,\ldots,r-1\}$. Then the map sending
$\alpha$ to $S(\alpha)$ is a bijection between the set of
compositions of $r$ and the set of subsets of $\{1,\ldots,r-1\}$.
For $S=\{i_1<\ldots<i_s\}\subset\{1,\ldots,r-1\}$, the inverse image
is given by $\alpha(S)=\{i_1,i_2-i_1,i_3-i_2,\ldots,r-i_s\}$. We
denote by $\C$ the set of compositions.

One may identify a composition with a {\it ribbon diagram}
\cite{MM}. For example, a composition $\alpha=(1,1,3,4,1,2)$ is
identified with
\begin{equation*}
\begin{array}{ccccccc}
    \bullet & & & & & & \\
    \bullet & & & & & & \\
    \bullet & \bullet & \bullet &  &  & & \\
    & &  \bullet & \bullet & \bullet & \bullet & \\
    &  &  &  & & \bullet & \\
    &  &  &  & & \bullet & \bullet \\
\end{array}\ \ \ .
\end{equation*}

A tableau $T$ obtained by filling a ribbon diagram $\alpha\in\C$
with entries in $\N$ is called a {\it quasi-ribbon tableau of shape
$\alpha$} if (1) the entries in each row (resp. column) are weakly
increasing from left to right (resp. from top to bottom), and (2)
the entries of positive integers (resp. half integers) are strictly
increasing in each column (resp. row). We denote by
$\B_{\N}(\alpha)$, or simply $\B(\alpha)$ if there is no confusion,
the set of quasi-ribbon tableaux of shape $\alpha$
(cf.~\cite{KT97}).

Suppose that $T\in \B(\alpha)$ is given. Let $w(T)$ be the word in
$\W$ obtained by reading the entries in $T$ row by row from top to
bottom, where in each row we read the entries from right to left.
For $i\in I$ and $x=e,f$, we define $x_iT$ to be the unique tableau
in $\B(\alpha)$ satisfying $w(x_i T)= x_i w(T)$, where we assume
that $x_i T={\bf 0}$ if $x_i w(T)={\bf 0}$.

Let $H_{\alpha}$ be the element in $\B(\alpha)$, which is defined in
the following way;
\begin{equation*}
\begin{array}{ccccccc}
    \hf & & & & & & \\
    \hf & & & & & & \\
    \hf & 1 & 1 &  & & &  \\
    & &  \frac{3}{2} & 2 & 2 & 2 & \\
    &  &  &  & & \frac{5}{2} & \\
    &  &  &  & & \frac{5}{2} & 3\\
\end{array}\ \ \ .
\end{equation*}
Then ${\rm wt}(H_{\alpha})\geq {\rm wt}(T)$ for all
$T\in\B(\alpha)$, where ${\rm wt}(T)={\rm wt}(w(T))$.

\begin{thm}\label{Bnalpha} For $\alpha\in\C$,
$\B(\alpha)$ is a $\frak{gl}_{\N}$-crystal and
$$\B(\alpha)=\{\,f_{i_1}\cdots f_{i_r} H_{\alpha}\,|\,r\geq 0,\ i_1,\ldots,i_r\in I \,\}\setminus\{{\bf 0}\}.$$
In particular, $\B(\alpha)$ is connected with a unique highest
weight element $H_{\alpha}$.
\end{thm}
\pf It is straightforward to check that $x_i T\in
\B(\alpha)\cup\{{\bf 0}\}$ is well-defined for $T\in\B(\alpha)$,
$x=e,f$, and $i\in I$. Hence, identifying $T$ with $w(T)$,
$\B(\alpha)$ becomes a subcrystal of $\W$.

Next, we will show that for each $T\in\B(\alpha)$, if $T\neq
H_{\alpha}$, then there exists $i\in I$ such that $e_i T\neq {\bf
0}$. Let $x$ be a node in $\alpha$, whose entry $r$ is not equal to
that of $H_{\alpha}$. We also assume that the other entries located
to the northwest of $x$, that is, the entries whose row or column
indices are no more than that of $x$, are equal to those in
$H_{\alpha}$. If $r\in\mathbb{Z}_{>0}$, then there is no $r-\hf$ to
the left of $x$ in the same row and no more $r$ in the rows strictly
lower than that of $x$. So we have $e_{r-\hf}T\neq {\bf 0}$.
Similarly, if $r\in\hf+\mathbb{Z}_{\geq 0}$, then there is no more
$r$ to the northwest of it and no $r-\hf$ to the right of $r$ in the
same row, which also implies that $e_{r-\hf}T\neq {\bf 0}$. This
completes the proof. \qed

\begin{figure}
\includegraphics{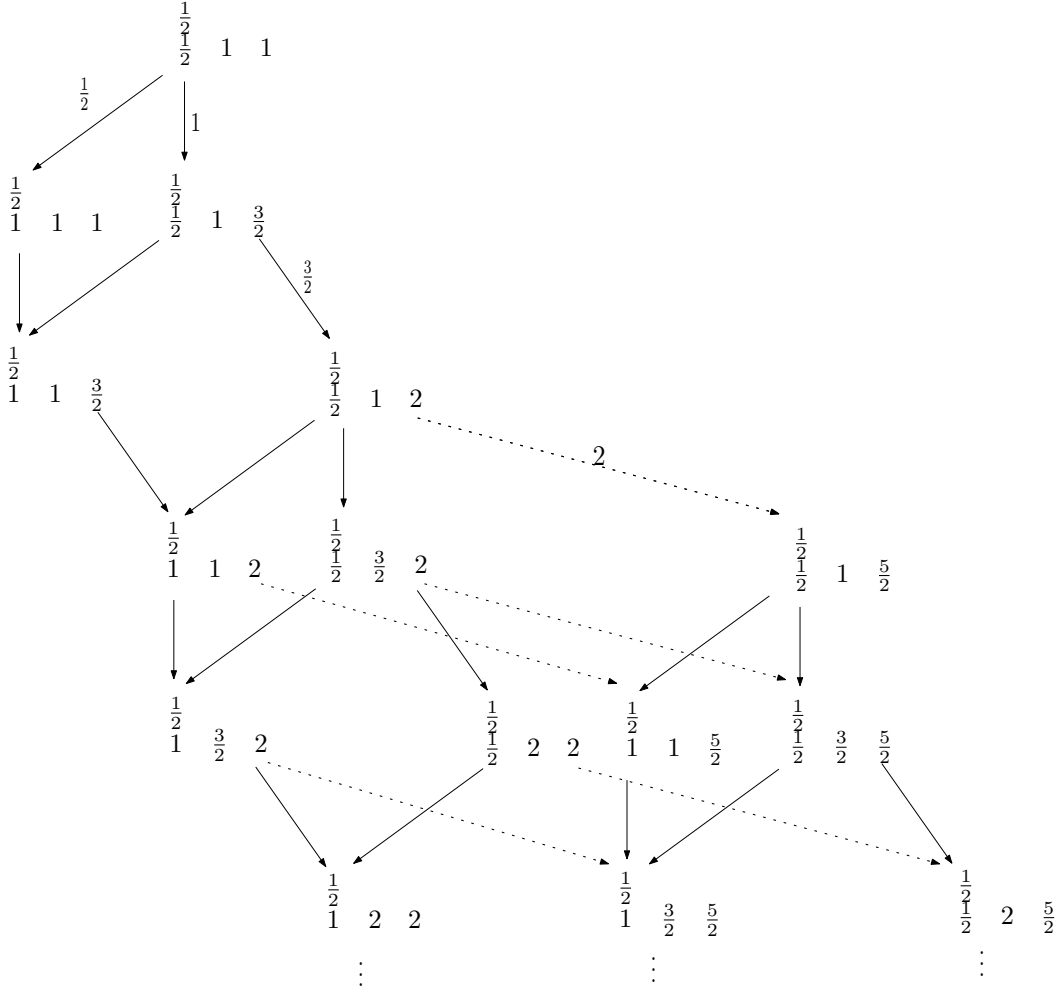}
\caption{$\B(1,3)$ with elements of height $\leq 6$ and arrows
colored by $\{\,\hf,1,\frac{3}{2},2\,\}$}
\end{figure}

\subsection{Crystal equivalence}\label{insertion}
Let $T\in\B(\alpha)$ be given for some  $\alpha\in\C$. For $k\geq
1$, let $t_k$ denote the $k$th entry of $T$, where we enumerate the
nodes in $T$ from northwest to southeast. Let $T^{\leq k}$ be the
sub-tableau consisting of the first $k$ nodes in $T$ and $T^{\geq
k}$ the sub-tableau obtained by removing $T^{\leq k-1}$ from $T$,
where $T^{\leq 0}$ is assumed to be the empty tableau.

Given $b\in \N$, choose the smallest $t_k$ such that $b\preccurlyeq
t_k$ (resp. $b\prec t_k$) if $b\in\mathbb{Z}_{>0}$ (resp.
$b\in\hf+\mathbb{Z}_{\geq 0}$). Then we define $b \rightarrow T$ to
be the quasi-ribbon tableau glueing $T^{\leq k-1}$, $T^{\geq k}$ and
$b$, where $b$ is placed below the last node of $T^{\leq k-1}$ and
to the left of the first node of $T^{\geq k}$ (cf.~\cite{KT97}).

\begin{ex}
\begin{equation*}
2\rightarrow
\begin{array}{cccc}
  1 & 1 & \frac{5}{2} &  \\
  & & \frac{5}{2} & 3
\end{array}
=
\begin{array}{cccc}
  1 & 1 &   &     \\
    & {\bf 2}  & \frac{5}{2} &     \\
  & & \frac{5}{2} &  3
\end{array},\quad \quad
\frac{5}{2}\rightarrow
\begin{array}{cccc}
  1 & 1 & \frac{5}{2} &  \\
  & & \frac{5}{2} & 3
\end{array}
=
\begin{array}{cccc}
  1 & 1 & \frac{5}{2}  &     \\
    &   & \frac{5}{2} &     \\
  & & {\bf \frac{5}{2}} &  3
\end{array}
\end{equation*}
\end{ex}

\begin{lem}\label{crystal equiv} For $T\in \B(\alpha)$ and $b\in\N$, we have $(b\rightarrow T)\simeq T\otimes b$.
\end{lem}
\pf It is enough to show that for $i\in I$ and $x=e,f$,
\begin{equation*}
x_i (b\rightarrow T)=
\begin{cases}
x_i b \rightarrow T, & \text{if $x_i (T\otimes b)=T \otimes (x_i
b)$},
\\
 b \rightarrow x_i T, & \text{if $x_i (T\otimes b)= (x_i T) \otimes b$},
\end{cases}
\end{equation*}
where we understand that ${\bf 0}\rightarrow T= b\rightarrow {\bf
0}={\bf 0}$.

We will prove the case when $x=f$ and $i,b\in\mathbb{Z}_{>0}$. The
other cases can be verified in a similar way. First, suppose that
$f_i (T\otimes b)=T \otimes f_i b$. Then $b=i$, and neither $i$ nor
$i+\hf$ appears as an entry of $T$. This implies that $f_i
(b\rightarrow T)=f_i b\rightarrow T$. Next, suppose that $f_i
(T\otimes b)=f_iT \otimes b$, and $t_k=i$, the $k$th entry of $T$
from northwest, becomes $i+\hf$ applying $f_i$ to $T$. If $b\neq
i,i+\hf$, then it is clear that $f_i (b\rightarrow T)= b\rightarrow
f_i T$. If $b=i$, then $b$ is located to the left of $t_k$  in the
same row of $b\rightarrow T$. If $b=i+\hf$, then $b$ is located
below $t_k$ in the same column of $b\rightarrow T$. This also
implies that $f_i (b\rightarrow T)= b\rightarrow f_i T$.  \qed

Now, for   $w=w_1\ldots w_r \in\W$, we define its {\it quasi
$P$-tableau} by
\begin{equation}
P(w)= w_r \rightarrow(\cdots\rightarrow(w_2 \rightarrow w_1)).
\end{equation}
By Theorem \ref{Bnalpha} and Lemma \ref{crystal equiv}, we have
\begin{cor}\label{component in Wn}
For $w\in\W$, we have $w\simeq P(w)$. In particular, each connected
component in $\W$ is isomorphic to $\B(\alpha)$ for some
$\alpha\in\C$.
\end{cor}

\begin{cor}\label{equivalence}
Each $w$ in $\W$ is $\frak{gl}_{\N}$-equivalent to a unique
quasi-ribbon tableau.
\end{cor}
\pf Let $T$ be a quasi-ribbon tableau which is
$\frak{gl}_{\N}$-equivalent to $w$. Then $T\simeq P(w)$, and hence
the highest weight vectors for the connected components of $T$ and
$P(w)$ are equal since they have the same weight. Since $T$ and
$P(w)$ are generated by the same highest weight vector, it follows
that $T=P(w)$. \qed

For $\alpha\in\C$ with $\alpha\vDash r$, a tableau $T$ obtained by
filling a ribbon diagram $\alpha$ with $\{1,\ldots,r\}$ is called a
{\it standard ribbon tableau of shape $\alpha$} if the entries in
each row are decreasing from left to right, and the entries in each
column are increasing from top to bottom.

For $w=w_1\ldots w_r \in\W$, we define its {\it quasi $Q$-tableau}
$Q(w)$ to be the standard ribbon tableau of the same shape as
$P(w)$, where we fill a node of $Q(w)$ with $i$ if the corresponding
position in $P(w)$ is given by $w_i$. Then as in the classical
Robinson-Schensted correspondence (cf.~\cite{Kn}), the map $w\mapsto
(P(w),Q(w))$ gives  a bijection from $\W$ to the set of pairs of a
quasi-ribbon tableau and a standard ribbon tableau of the same shape
(cf.~\cite{KT97}). Furthermore, it is straightforward to check that
$Q(x_iw)=Q(w)$ whenever $x_iw\neq {\bf 0}$ for $i\in I$ and $x=e,f$
(cf.~ Proposition 4.17 in \cite{KK}, and \cite{LT}).
\begin{ex}{\rm Let $w=1\ \hf \ 1\ \frac{5}{2}\ 2\ 2$. Then
\begin{equation*}
P(w)=
\begin{array}{cccc}
  \hf &  &  &  \\
  1   & 1 &  &  \\
   & 2 & 2 & \frac{5}{2}
\end{array}
,\ \ \ \ Q(w)=
\begin{array}{cccc}
  2 &  &  &  \\
  3   & 1 &  &  \\
   & 6 & 5 & 4
\end{array}.
\end{equation*}}
\end{ex}
Summarizing the above arguments, we have
\begin{thm}
For $\alpha\in\C$, let $RT(\alpha)$ be the set of standard ribbon
tableaux of shape $\alpha$. For $T\in RT(\alpha)$,
$\B(T)=\{\,w\in\W\,|\,Q(w)=T\,\}$ is isomorphic to $\B(\alpha)$.
Hence, as a $\frak{gl}_{\N}$-crystal, we have
\begin{equation*}
\W = \bigoplus_{\alpha\in\C}\bigoplus_{T\in RT(\alpha)} \B(T).
\end{equation*}
\end{thm}

\subsection{Stability of crystal graphs}\label{stablity of crystals}
For $n\in\mathbb{Z}_{>0}$, put $\N^{\preccurlyeq
n}=\{\,r\in\N\,|\,r\preccurlyeq n\,\}$.  Then $\B_{\N^{\preccurlyeq
n}}(2^{n-1},1)$ has $2^{2n-1}$ elements given by
\begin{equation}\label{elementsinzigzag}
f_{1/2}^{m_{1/2}}f_{1}^{m_{1}}\cdots f_{n-\hf}^{m_{n-\hf}}
H_{(2^{n-1},1)},
\end{equation}
where $m_i=0,1$ for $\alpha_i\in \Delta_{\N^{\preccurlyeq n}}$. Note
that each element in $\B_{\N^{\preccurlyeq n}}(2^{n-1},1)$ is
uniquely determined by a sequence $(m_{\hf},\ldots,m_{n-\hf})\in
\{0,1\}^{2n-1}$ (see Figure 2).

\begin{figure}
\includegraphics{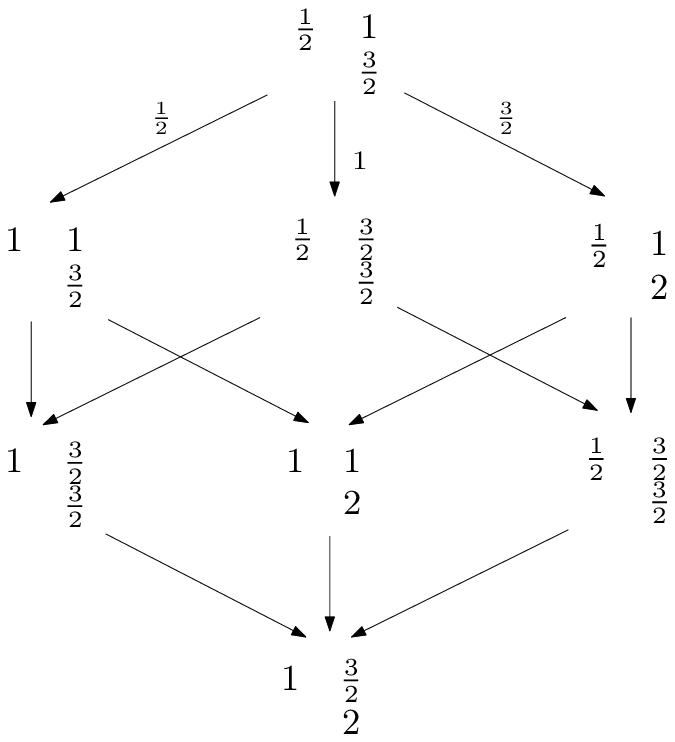}
\caption{Crystal graph of $\B_{\N^{\preccurlyeq 2}}(2,1)$}
\end{figure}

Let us say that a node in a ribbon diagram is a {\it corner} if it
is the leftmost or rightmost node in each row having at least two
nodes, or it is the last node in the diagram when enumerated from
northwest to southeast. Note that for $\alpha\in\C$,
$\B_{\N^{\preccurlyeq n}}(\alpha)$ is non-empty if and only if
$\alpha$ has at most $2n$ corners. Suppose that $\alpha$ has either
$2n-1$ or $2n$ corners, equivalently, $({\rm
wt}(H_{\alpha}),\epsilon_r)\neq 0$ for $r\prec n$, and $0$ for
$r\succ n$. Then each $T$ in $\B_{\N^{\preccurlyeq n}}(\alpha)$
differs from $H_{\alpha}$ only in corners. Hence, it is not
difficult to see that $\B_{\N^{\preccurlyeq n}}(\alpha)$ consists of
$2^{2n-1}$ elements given by $f_{1/2}^{m_{1/2}}f_{1}^{m_{1}}\cdots
f_{n-\hf}^{m_{n-\hf}} H_{\alpha}$, where $m_i=0,1$ for $\alpha_i\in
\Delta_{\N^{\preccurlyeq n}}$. This implies the following fact.

\begin{prop}\label{stablity} Let $\alpha$ be a composition
whose number of corners is either $2n-1$ or $2n$. Then there exists
a unique bijection
$$\theta : \B_{\N^{\preccurlyeq n}}(2^{n-1},1) \longrightarrow \B_{\N^{\preccurlyeq n}}(\alpha),$$
which maps $H_{(2^{n-1},1)}$ to $H_{\alpha}$ and commutes with $x_i$
for $x=e,f$ and $\alpha_i\in \Delta_{\N^{\preccurlyeq n}}$. Here we
assume that $\theta({\bf 0})={\bf 0}$. In other words,
$\B_{\N^{\preccurlyeq n}}(2^{n-1},1)$ and $ \B_{\N^{\preccurlyeq
n}}(\alpha)$ are isomorphic as
$\Delta_{\N^{ \preccurlyeq n}}$-colored oriented graphs. 
\end{prop}

\section{Enriched $P$-partitions}

\subsection{Crystal structure on the set of enriched $P$-partitions}
Let us recall the notion of an enriched $(P,\gamma)$-partition
introduced by Stembridge \cite{St97}. We follow the notations in
\cite{St97} with a little modification. Let $P=(X,<)$ be a finite
set $X$ with a partial ordering $<$. Let $\gamma : X \rightarrow
\mathbb{Z}_{>0}$ be an injective map, which will be called a
labeling of $X$. We call a pair $(P,\gamma)$ a labeled poset. Then
an {\it enriched $(P,\gamma)$-partition} is a map $\sigma :
X\rightarrow \N$ such that for all $x<y$ in $P$,
\begin{itemize}
\item[(1)] $\sigma(x)\preccurlyeq\sigma(y)$,

\item[(2)] $\sigma(x)=\sigma(y)$ and $\sigma(x)\in\mathbb{Z}_{>0}$ implies
$\gamma(x)<\gamma(y)$,

\item[(3)] $\sigma(x)=\sigma(y)$ and $\sigma(x)\in\hf+\mathbb{Z}_{\geq 0}$ implies
$\gamma(x)>\gamma(y)$.
\end{itemize}
We denote by $\mathcal{E}(P,\gamma)$ the set of enriched
$(P,\gamma)$-partitions.

Suppose that $|X|=r$. Define an embedding $\psi :
\mathcal{E}(P,\gamma) \rightarrow \W$ by
$\psi(\sigma)=\sigma(x_1)\cdots \sigma(x_r)$ for $\sigma \in
\mathcal{E}(P,\gamma)$, where $x_i\in X$ ($1\leq i\leq r$) are
arranged so that $\gamma(x_1)>\gamma(x_2)>\cdots>\gamma(x_r)$.

\begin{thm} The image of $\mathcal{E}(P,\gamma)$ under $\psi$ together with $\{{\bf
0}\}$ is stable under $e_i,f_i$ for $i\in I$. Hence,
$\mathcal{E}(P,\gamma)$ becomes a $\frak{gl}_{\N}$-crystal.
\end{thm}
\pf We will prove that $x_i\psi(\sigma)\in
\psi(\mathcal{E}(P,\gamma))$ for $x=e,f$ and $i\in I$, when
$x_i\psi(\sigma)\neq {\bf 0}$. We assume that $x=f$ and
$i\in\mathbb{Z}_{>0}$ since the other cases can be checked in the
same way.

Suppose that
\begin{equation*}
f_i \psi(\sigma)=\sigma(x_1)\cdots  f_i \sigma(x_s)\cdots
\sigma(x_r)\neq {\bf 0},
\end{equation*}
where $\sigma(x_s)=i$ and hence $f_i \sigma(x_s)=i+\hf$. Recall that
$s$ is the smallest index such that $(\alpha_i,{\rm
wt}(\sigma(x_s)))\neq 0$. Let $\tau : X \rightarrow \N$ be defined
by $\tau(x_t)=\sigma(x_t)$ for $t\neq s$ and $\tau(x_s)=i+\hf$.

Suppose that $x_s<x_t$ in $P$. If $\sigma(x_s)=\sigma(x_t)=i$, then
$\gamma(x_s)<\gamma(x_t)$, which contradicts the minimality of $s$.
If $\sigma(x_s)=i\prec\sigma(x_t)=i+\hf$, then we have $\tau(x_s)=
\tau(x_t)$, and $\gamma(x_s)>\gamma(x_t)$ by the minimality of $s$.
If $\sigma(x_s)=i\prec i+\hf\prec\sigma(x_t)$, then we have
$\tau(x_s)\prec \tau(x_t)$. Next, suppose that $x_t<x_s$ in $P$.
Then it is clear that $\tau(x_t)\prec \tau(x_s)$. Hence, it follows
that $\tau\in \mathcal{E}(P,\gamma)$ and
$\psi(\tau)=f_i\psi(\sigma)$. \qed

We mean a linear extension of $P$ by a total ordering
$w=\{\,w_1<\cdots<w_r\,\}$ on $P$ preserving its partial ordering.
We denote by $\mathcal{L}(P)$ the set of linear extensions of $P$.
For $w\in \mathcal{L}(P)$, let
$D(w,\gamma)=\{\,i\,|\,\gamma(w_i)>\gamma(w_{i+1})\,\}$ be the
descent of $w$ with respect to $\gamma$. Then we have
\begin{equation}
\begin{split}
\mathcal{E}(w,\gamma)=\{\,\sigma : X \rightarrow \N\,|\,& (1)\
\sigma(w_1)\preccurlyeq\cdots\preccurlyeq \sigma(w_r), \\
& (2)\ \sigma(w_i)=\sigma(w_{i+1})\in\mathbb{Z}_{>0} \Rightarrow
i\not\in
D(w,\gamma), \\
& (3)\ \sigma(w_i)=\sigma(w_{i+1})\in\hf+\mathbb{Z}_{\geq 0}
\Rightarrow i \in D(w,\gamma)\,\}.
\end{split}
\end{equation}

\begin{lem} For $w\in \mathcal{L}(P)$, $\mathcal{E}(w,\gamma)\simeq\B(\alpha)$ as
a $\frak{gl}_{\N}$-crystal, where $\alpha=\alpha(D(w,\gamma))\in\C$.
\end{lem}
\pf Let us identify $w_k$ ($1\leq k\leq r$) with the $k$th node in
the ribbon diagram $\alpha$ from its northwest. This induces an
isomorphism of $\frak{gl}_{\N}$-crystals from
$\mathcal{E}(w,\gamma)$ to $\B(\alpha)$. \qed

\begin{rem}{\rm
Note that the crystal graph structure on $\mathcal{E}(w,\gamma)$
depends only on $D(w,\gamma)$. }
\end{rem}

Let us consider the decomposition of $\mathcal{E}(P,\gamma)$ as a
$\frak{gl}_{\N}$-crystal. Given $\sigma\in \mathcal{E}(P,\gamma)$,
we define a linear extension $w$ on $P$ as follows;
\begin{itemize}
\item[(1)] arrange the elements of $X$ in increasing order of their
values of $\sigma$,

\item[(2)] if there exist elements $x$ in $X$ with the same value
$\sigma(x)\in\mathbb{Z}_{>0}$ (resp.
$\sigma(x)\in\hf+\mathbb{Z}_{\geq 0}$), then arrange them in order
of their increasing (resp. decreasing) values of $\gamma$.
\end{itemize}
By definition, we have $\sigma\in \mathcal{E}(w,\gamma)$, and this
correspondence induces a bijection $\pi : \mathcal{E}(P,\gamma)
\rightarrow \bigsqcup_{w\in\mathcal{L}(P)}\mathcal{E}(w,\gamma)$
(Lemma 2.1 in \cite{St97}). Furthermore, it is straightforward to
see that $\pi$ commutes with $x_i$ for $x=e,f$ and $i\in I$, where
we assume that $\pi({\bf 0})={\bf 0}$ and $x_i{\bf 0}={\bf 0}$.
Hence, we obtain the following.

\begin{cor}\label{FLEPP} Given a labeled poset $(P,\gamma)$, we have
\begin{equation*}
\mathcal{E}(P,\gamma) \simeq
\bigoplus_{w\in\mathcal{L}(P)}\mathcal{E}(w,\gamma),
\end{equation*}
as a $\frak{gl}_{\N}$-crystal. That is, for each $\alpha\in\C$, the
multiplicity of $\B(\alpha)$ in $\mathcal{E}(P,\gamma)$ is equal to
the number of $w\in \mathcal{L}(P)$ such that
$\alpha(D(w,\gamma))=\alpha$. 
\end{cor}

\subsection{Decomposition of $\frak{gl}_{\N}$-crystals}

Consider $\B(\lambda)=\B_{\N}(\lambda)$ for $\lambda\in\cP$ (see
Section 2.2). We assume that $\lambda$ is equipped with a partial
ordering  given by $\lambda(i,j)<\lambda(i,j+1)$ and
$\lambda(i,j)<\lambda(i+1,j)$, where $\lambda(i,j)$ denotes the node
in $\lambda$ located in the $i$th row from the top and the $j$th
column from the left. We consider a labeled poset
$(\lambda,\gamma)$, where $\gamma$ satisfies
$\gamma(\lambda(i,j))<\gamma(\lambda(i,j+1))$ and
$\gamma(\lambda(i,j))>\gamma(\lambda(i+1,j))$. Then
$\mathcal{E}(\lambda,\gamma)$ is equal to $\B(\lambda)$ as a set,
and if we use the admissible reading on $\B(\lambda)$ given by the
reverse ordering of $\gamma$-values, then the actions of $e_i,f_i$
($i\in I$) on both sets coincide.

Now, let $ST(\lambda)$ be the set of standard tableaux of shape
$\lambda$, that is, the set of order preserving bijections $T:
\lambda \rightarrow \{\,1,\ldots,r\,\}$, where we regard $\lambda$
as a poset with $r$ elements. Then $\mathcal{L}(\lambda)$ can be
identified with $ST(\lambda)$. For $T\in ST(\lambda)$, we denote by
$D(T)$ the descent set of $T$ as a linear extension of $\lambda$.
Hence, we have $k\in D(T)$ if and only if the column index of
$T^{-1}(k)$ is greater than or equal to that of $T^{-1}(k+1)$. We
put $\alpha(T)=\alpha(D(T))$ for simplicity. Hence, by Corollary
\ref{FLEPP}, we obtain the following decomposition of $\B(\lambda)$.

\begin{prop}\label{decomposition of Blambda} For $\lambda\in\cP$, we have
\begin{equation*}
\B(\lambda) \simeq \bigoplus_{T\in ST(\lambda)}\B(\alpha(T)).
\end{equation*}
\end{prop}

\begin{rem}\label{decomposition of Blambda skew version}{\rm
Proposition \ref{decomposition of Blambda} can be directly extended
to the case of $\B(\lambda/\mu)$ for a skew Young diagram
$\lambda/\mu$. Moreover, for a strict partition $\lambda$, the set
of shifted $\N$-tableaux of shape $\lambda$ is also a
$\frak{gl}_{\N}$-crystal, and we have a similar decomposition (cf.
\cite{St97}). }
\end{rem}

Let $\alpha\in\C$ be given with $\alpha\vDash r$. Let
$x_1,\ldots,x_r$ denote the nodes in the diagram of $\alpha$
enumerated from northwest to southeast. We assume that $\alpha$ is
equipped with a total ordering  $w_{\alpha}=\{\,x_1<\cdots<x_r\,\}$.
We also identify $T\in \B(\alpha)$ with a function $T : \alpha
\rightarrow \N$. A canonical labeling $\gamma_{\alpha}$ of $\alpha$
is defined by $\gamma_{\alpha}(x_i)=r-k+1$ if $T(x_i)$ is read in
the $k$th letter of $\psi(T)\in \W$ for $T\in\B(\alpha)$. Then
clearly we have $\B(\alpha)=\mathcal{E}(w_{\alpha},\gamma_{\alpha})$
and the actions of $x_i$ for $x=e,f$ and $i\in I$ on both sets
coincide. For $t\geq 0$, we define $\gamma_{\alpha}^{[t]}$ by
$\gamma_{\alpha}^{[t]}(x_i)=\gamma_{\alpha}(x_i)+t$ for $1\leq i\leq
r$.

Suppose that $\alpha,\beta\in\C$ are given with $\beta\vDash s$.
Consider a labeled poset $(w_{\alpha}\cup
w_{\beta},\gamma_{\alpha}^{[s]}\cup \gamma_{\beta})$, which is a
disjoint union of labeled posets. Then we can check that
$\B(\alpha)\otimes \B(\beta)$ is isomorphic to
$\mathcal{E}(w_{\alpha}\cup w_{\beta},\gamma_{\alpha}^{[s]}\cup
\gamma_{\beta})$ as  a $\frak{gl}_{\N}$-crystal. By Corollary
\ref{FLEPP}, we obtain the following tensor product decomposition.

\begin{prop}\label{LR rule} For $\alpha,\beta\in\C$ with
$\beta\vDash s$, we have
\begin{equation*}
\B(\alpha)\otimes \B(\beta)\simeq \bigoplus_{w\in
\mathcal{L}(w_{\alpha}\cup
w_{\beta})}\mathcal{E}(w,\gamma_{\alpha}^{[s]}\cup \gamma_{\beta}).
\end{equation*}
\end{prop}
Note that a linear extension of $w_{\alpha}\cup w_{\beta}$ is called
a shuffle of $\alpha$ and $\beta$.

\subsection{RSK correspondence}
Let
\begin{equation}
\begin{split}
\Omega=\{\,& (\bi,\bj)\in
\W\times \W\,| \\
& \text{(1) $\bi=i_1\cdots i_r$ and $\bj=j_1\cdots j_r$  for some $r\geq 0$}, \\
& \text{(2) $(i_1,j_1)\leq \cdots \leq (i_r,j_r)$}, \\
& \text{(3) $i_k-j_k\not\in \mathbb{Z}$ implies $(i_k,j_k)\neq
(i_{k\pm1},j_{k\pm 1})$}\},
\end{split}
\end{equation}
where for $(i,j)$ and $(k,l)\in \N\times \N$, the {\it super
lexicographic ordering} is given by
\begin{equation}\label{partial order}
(i,j)< (k,l) \ \ \ \ \Longleftrightarrow \ \ \ \
\begin{cases}
(j<l) & \text{or}, \\
(j=l\in\mathbb{Z}_{>0},\ \text{and} \ i>k) & \text{or}, \\
(j=l\in\hf+\mathbb{Z}_{\geq 0},\ \text{and} \ i<k) &.
\end{cases}
\end{equation}
Also, let $\Omega^*$ be the set of pairs $(\bk,\bl)\in \W\times \W$
such that $(\bl,\bk)\in \Omega$.

Given $(\bi,\bj)\in \Omega$, we define $x_i (\bi,\bj)=(x_i\bi,\bj)$
for $x=e,f$ and $i\in I$, where we assume that $x_i (\bi,\bj)={\bf
0}$ if $x_i\bi={\bf 0}$, and set ${\rm wt}(\bi,\bj)={\rm wt}(\bi)$.
Similarly, given $(\bk,\bl)\in \Omega^*$, we define $x_j^*
(\bk,\bl)=(\bk,x_j\bl)$ for $x=e,f$ and $j\in I$, and set ${\rm
wt}^*(\bk,\bl)={\rm wt}(\bl)$. Then as in \cite{K07}, we can check
that

\begin{lem}[cf.~\cite{K07} Lemma 3.1]\label{crystal Omega} Under the above hypothesis,
$\Omega$ and $\Omega^*$ are $\frak{gl}_{\N}$-crystals with respect
to $x_i$ and $x_i^*$ for $x=e,f$ and $i\in I$, respectively. 
\end{lem}

Consider
\begin{equation}
\begin{split}
\M= \{\,A=&(a_{rs})_{r,s\in \N}\,|\, \\
& \text{(1) $a_{rs}=0$ for all sufficiently large $r$ and $s$, } \\
& \text{(2) $a_{rs}\in\mathbb{Z}_{\geq 0}$,\ and \   $a_{rs}\leq 1$
unless $r-s\in \mathbb{Z}$}\, \}.
\end{split}
\end{equation}

For $(\bi,\bj)\in \Omega$, define $A(\bi,\bj)=(a_{rs})$ to be the
matrix in $\M$, where $a_{rs}$ is the number of $k$'s such that
$(i_k,j_k)=(r,s)$ for $r,s\in \N$. Then, it follows that the map
$(\bi,\bj)\mapsto A(\bi,\bj)$ gives a bijection from $\Omega$ to
$\M$, where the pair of empty words $(\emptyset,\emptyset)$
corresponds to zero matrix. Similarly, we have a bijection
$(\bk,\bl)\mapsto A(\bk,\bl)$ from $\Omega^*$ to $\M$. With these
bijections, $\M$ becomes a $\frak{gl}_{\N}$-crystal with respect to
both $x_i$ and $x_i^*$ for $x=e,f$ and $i\in I$ by Lemma
\ref{crystal Omega}. For convenience, let us say that $\M$ is a
$\frak{gl}^*_{\N}$-crystal when we consider its crystal structure
with respect to $x_i^*$.

\begin{lem}[cf.~\cite{K07} Lemma 3.4]\label{bicrystal}
$\M$ is a $(\frak{gl}_{\N},\frak{gl}^*_{\N})$-bicrystal, that is,
$e_i,f_i$ commute with $e_j^*,f_j^*$ for $i,j\in I$, where we assume
that $x_i{\bf 0}=x_j^*{\bf 0}={\bf 0}$ for $x=e,f$. 
\end{lem}

Given $A\in \M$, suppose that $A=A(\bi,\bj)=A(\bk,\bl)$ for
$(\bi,\bj)\in \Omega$ and $(\bk,\bl)\in \Omega^*$. We define
\begin{equation}\label{varpi}
\varpi(A)=(P_1(A),P_2(A))=({P}(\bi),{P}(\bl)).
\end{equation}
Note that $A$ is $\frak{gl}_{\N}$ (resp.
$\frak{gl}^*_{\N}$)-equivalent to $P_1(A)$ (resp. $P_2(A)$).
\begin{prop}\label{RSK} The map $\varpi$ induces an isomorphism of
$(\frak{gl}_{\N},\frak{gl}^*_{\N})$-bicrystals
\begin{equation*}
\varpi : \M \longrightarrow
\bigsqcup_{\lambda\in\cP}\bigsqcup_{P,Q\in
ST(\lambda)}\B(\alpha(P))\times \B(\alpha(Q)).
\end{equation*}
\end{prop}
\pf Given $A\in \M$, suppose that $A=A(\bi,\bj)=A(\bk,\bl)$ for
$(\bi,\bj)\in \Omega$ and $(\bk,\bl)\in \Omega^*$. Define
$\pi(A)=({\bf P}(\bi),{\bf P}(\bl))$ (see \eqref{column insertion}).
One can extend the arguments for Young tableaux (see, for example,
$\S$4.2 in \cite{Fu}) to the super case without difficulty to prove
that $\pi$ induces a bijection from $\M$ to $\bigsqcup_{\lambda\in
\cP} \B(\lambda)\times \B(\lambda)$ (cf.~\cite{Kn}).

Let $\omega$ be the isomorphism given in Proposition
\ref{decomposition of Blambda}. Then we have a bijection $\pi' : \M
\rightarrow \bigsqcup_{\lambda\in\cP }\bigsqcup_{P,Q\in
ST(\lambda)}\B(\alpha(P))\times \B(\alpha(Q))$ defined by sending
$A$ to $(\omega({\bf P}(\bi)),\omega({\bf P}(\bl)))$. Furthermore,
$\omega({\bf P}(\bi))$ and $\omega({\bf P}(\bl))$ are quasi-ribbon
tableaux equivalent to $\bi$ and $\bl$ respectively, and by
Corollary \ref{equivalence}, it follows that $\omega({\bf
P}(\bi))=P(\bi)$, $\omega({\bf P}(\bl))=P(\bl)$ and hence
$\varpi=\pi'$.

Suppose that $A\in\M$ is given. If $x_j^*A\neq {\bf 0}$ for some
$x=e,f$ and $j\in I$, then $A$ is $\frak{gl}_{\N}$-equivalent to
$x_j^*A$ (cf.~\cite{K07} Lemma 3.5), and by Lemma \ref{equivalence},
we have $P_1(x_j^*A)=P_1(A)$. Similarly, if $x_iA\neq {\bf 0}$ for
some $x=e,f$ and  $i\in I$, then $P_2(x_iA)=P_2(A)$. This implies
that $\varpi$ is a morphism of
$(\frak{gl}_{\N},\frak{gl}^*_{\N})$-bicrystals. \qed

Let $\mathfrak{S}_k$ be the symmetric group on $k$ letters. For
$\sigma\in \mathfrak{S}_k$, let
$D(\sigma)=\{\,i\,|\,\sigma(i)>\sigma(i+1)\,\}$ be the descent of
$w$. Put $\alpha(\sigma)=\alpha(D(\sigma))$. Let $(P,Q)$ be the pair
of standard tableaux of the same shape, which corresponds to
$\sigma$ under the classical Robinson-Schensted correspondence. Then
$D(P)=D(\sigma)$ and $D(Q)=D(\sigma^{-1})$ \cite{Sch}, that is,
$\alpha(P)=\alpha(\sigma)$ and $\alpha(Q)=\alpha(\sigma^{-1})$.

Combining with Proposition \ref{RSK}, we obtain another version of
the Gessel's result \cite{G} in terms of crystal graphs.

\begin{thm}[cf.~\cite{G}]\label{Gessel} Let $\M^k=\{\,A=(a_{rs})\in\M\,|\,\sum_{r,s}a_{rs}=k\,\}$ for
$k\in\mathbb{Z}_{>0}$. Then $\varpi$ restricts to the following
isomorphism of $(\frak{gl}_{\N},\frak{gl}^*_{\N})$-bicrystals;
\begin{equation*}
\varpi : \M^k \longrightarrow
\bigsqcup_{\sigma\in\frak{S}_k}\B(\alpha(\sigma))\times
\B(\alpha(\sigma^{-1})).
\end{equation*}
\end{thm}

As a corollary, we have an interesting application to permutation
enumeration.
\begin{cor}[cf.~\cite{G}]
Given $S,S'\subset \{1,\ldots,k-1\}$, the number of permutations
$\sigma\in\frak{S}_k$ satisfying $D(\sigma)=S$ and
$D(\sigma^{-1})=S'$ is equal to the number of matrices $A\in \M^k$
satisfying $e_i A=e^*_jA={\bf 0}$ for all $i,j\in I$ with ${\rm
wt}(A)={\rm wt}(H_{\alpha(S)})$ and ${\rm wt}^*(A)={\rm
wt}(H_{\alpha(S')})$.
\end{cor}

\section{Non-standard Borel subalgebras and branching rule}
For $m\in\mathbb{Z}_{\geq 0}$, let us consider the crystal graphs
for $\frak{gl}_{\mathcal{N}(m)}$, where
$$\mathcal{N}(m)=\Big{\{}\,-m\prec   \cdots \prec -1 \prec\ \text{\tiny $\hf$} \prec 1\prec \text{\tiny $\frac{3}{2}$} \prec
2 \prec \cdots\,\Big{\}}.$$ As usual,
$\mathcal{N}(m)_0=\mathcal{N}(m)\cap\mathbb{Z}$ and
$\mathcal{N}(m)_1=\mathcal{N}(m)\cap \left(\hf+\mathbb{Z}\right)$.
We assume that $\N(0)=\N$. Then
$\Delta_{\mathcal{N}(m)}=\{\,\alpha_{i}=\epsilon_{i-1}-\epsilon_{i}\
(-m+1\leq i\leq -1),\ \alpha_0=\epsilon_{-1}-\epsilon_{\hf},\
\alpha_r=\epsilon_r-\epsilon_{r+\hf}\ (r\in\hf\mathbb{Z}_{>0})\,\}$.
Denote by $I(m)$ the index set for simple roots. The associated
Dynkin diagram is given by \vskip 5mm
\begin{center}
\includegraphics{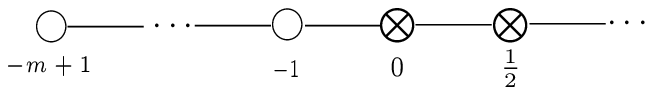}.
\end{center}
For $n\in\N(m)$, we put $\N(m)^{\preccurlyeq
n}=\{\,r\in\N(m)\,|\,r\preccurlyeq n\,\}$.

\subsection{Highest weight crystals}
Let $\C(m)$ be the set of pairs $(\lambda,\alpha)$ where
$\lambda=(\lambda_1,\ldots,\lambda_m)$ is a partition with length at
most $m$, and $\alpha=(\alpha_1,\ldots,\alpha_r)$ is a composition
such that $\alpha$ is non-empty only if $\lambda_m\neq 0$. We assume
that $\C(0)=\C$. One may identify $(\lambda,\alpha)\in \C(m)$ with a
diagram obtained by placing the first node of $\alpha$ from
northwest right below the leftmost node in the $m$th row of
$\lambda$. For example, when $m=3$, $\lambda=(5,4,2)$, and
$\alpha=(1,3,4)$, the corresponding diagram is
\begin{equation*}
\begin{array}{ccccccc}
  \bullet & \bullet & \bullet &  \bullet  & \bullet \\
  \bullet & \bullet & \bullet &  \bullet  &\\
  \bullet & \bullet &  &  &\\
  \bullet &  & &  & \\
  \bullet & \bullet & \bullet & &\\
  &  & \bullet & \bullet & \bullet & \bullet
\end{array}\ \ .
\end{equation*}
We call $\lambda$ the {\it body}, and $\alpha$ the {\it tail} of the
diagram of $(\lambda,\alpha)$. The first node of the tail from
northwest will be called the {\it joint} of $(\lambda,\alpha)$.

Let $\B_{\mathcal{N}(m)}(\lambda,\alpha)=\B(\lambda,\alpha)$ be the
set of tableaux $T$ obtained by filling the diagram
$(\lambda,\alpha)\in\C(m)$ with entries in $\mathcal{N}(m)$ such
that
\begin{itemize}
\item[(1)] $T$ is semistandard in the usual sense,

\item[(2)] if $b$ is the entry of its joint and $b\in\mathbb{Z}_{>0}$, then all the
entries in the body are smaller than $b$.

\item[(3)] if $b$ is the entry of its joint and $b\in\hf+\mathbb{Z}_{\geq 0}$, then all the
entries in the body are smaller than or equal to $b$.
\end{itemize}
We call $T\in \B(\lambda,\alpha)$ a {\it semistandard tableaux of
shape $(\lambda,\alpha)$}. We define $H_{(\lambda,\alpha)}$ to be
the tableau obtained by gluing $H_{\lambda}$ and $H_{\alpha}$. For
$T\in \B(\lambda,\alpha)$, let $w(T)$ be the word in
$\W_{\mathcal{N}(m)}=\W$ obtained from $T$ with respect to row
reading. Then for $i\in I(m)$ and $x=e,f$, we define $x_i T$ to be
the tableau of shape $(\lambda,\alpha)$ corresponding to $x_iw(T)$.

\begin{prop}\label{Blambdaalpha} For $(\lambda,\alpha)\in\C(m)$, $\B(\lambda,\alpha)$ is a
$\frak{gl}_{\mathcal{N}(m)}$-crystal, and
\begin{equation*}
\B(\lambda,\alpha)=\{\,f_{i_1}\cdots f_{i_r}
H_{(\lambda,\alpha)}\,|\,r\geq 0,\ i_1,\ldots,i_r\in I(m)
\,\}\setminus\{{\bf 0}\}.
\end{equation*}
In particular, $\B(\lambda,\alpha)$ is connected with a unique
highest weight element $H_{(\lambda,\alpha)}$.
\end{prop}
\pf We can check that $x_i T\in \B(\lambda,\alpha)\cup\{{\bf 0}\}$
for $x=e,f$, $i\in I(m)$, and $T\in \B(\lambda,\alpha)$. Hence
$\B(\lambda,\alpha)$ is a $\frak{gl}_{\mathcal{N}(m)}$-subcrystal of
$\W$ with respect to row reading. We leave the details to the
readers.

Next we claim that for $T\in \B(\lambda,\alpha)$, if $e_iT={\bf 0}$
for all $i\in I(m)$, then $T=H_{(\lambda,\alpha)}$. Let $k$ be the
number of positive entries lying in the body of $T$. If $k=0$, it is
clear that $T=H_{(\lambda,\alpha)}$ since the body and the tail of
$T$ should be $H_{\lambda}$ and $H_{\alpha}$ by Theorem
\ref{Bmnlambda} and \ref{Bnalpha}, respectively. Suppose that $k>0$.
Since $e_iT={\bf 0}$ for all $i\geq 1$, there exists at least one
$\hf$ in $T$. Also, by the condition (3), there exists at least one
$\hf$ in the body of $T$. Choose the $\hf$ which is located in the
highest position in the body of $T$. Note that no $-1$ can be read
before such $\hf$ with respect to row reading of $T$ since
$e_iT={\bf 0}$ for $-m+1\leq i\leq -1$ and hence all $-1$ should lie
in the $m$th row of $T$. Then applying $e_0$ to $T$ changes such
$\hf$ to $-1$, which is a contradiction. This completes the proof.
\qed

\begin{cor}\label{equivalence for nonstandard}
Each $w$ in $\W$ is $\frak{gl}_{\mathcal{N}(m)}$-equivalent to a
unique semistandard tableau of shape $(\lambda,\alpha)\in\C(m)$.
\end{cor}
\pf It suffices to show that each $w$ in $\W$ is
$\frak{gl}_{\mathcal{N}(m)}$-equivalent to a tableau $T\in
\B(\lambda,\alpha)$ for some $(\lambda,\alpha)\in\C(m)$. Then the
uniqueness follows from Proposition \ref{Blambdaalpha} together with
the same arguments as in Corollary \ref{equivalence}.

First, consider $T={\bf P}(w)$ which is
$\frak{gl}_{\mathcal{N}(m)}$-equivalent to $w$ (cf.\eqref{column
insertion}). Let $b$ be the entry of $T$ located in the first column
of the $(m+1)$st row, which is obviously greater than $-1$. Let
$T_2$ be the subtableau of $T$ consisting of entries $b'$ such that

\begin{itemize}
\item[(1)] $b\preccurlyeq b'$,

\item[(2)] $b=b'$ only when $b\in\mathbb{Z}_{>0}$ or $b'$ is
located below $b$.
\end{itemize}
Let $T_1$ be the complement of $T_2$ in $T$. By Corollary
\ref{equivalence}, $T_2$ is equivalent to a unique quasi-ribbon
tableau $T_2'$ as a $\frak{gl}_{\N}$-crystal element.

Let ${\bf P}'(w)$ be the tableau obtained by placing the first node
of $T_2'$ from northwest right below the leftmost node in the $m$th
row of $T_1$. By construction, we have ${\bf P}'(w)\in
\B(\lambda,\alpha)$ for some $(\lambda,\alpha)\in\C(m)$. Clearly,
$w$ is $\frak{gl}_{[m|0]}$-equivalent to ${\bf P}'(w)$. Let $T_1'$
be the subtableau of $T_1$ consisting of positive entries, then $w$
is $\frak{gl}_{\N}$-equivalent to $T_1'\otimes T_2'$ and hence to
${\bf P}'(w)$ by Proposition \ref{decomposition of Blambda}. In
particular, we have ${\bf P}'(x_i w)=x_i{\bf P}'(w)$ for $x=e,f$ and
$i\neq 0$. Finally, we can check that $x_0{\bf P}(w)\neq {\bf 0}$ if
and only if $x_0{\bf P}'(w)\neq {\bf 0}$, and $x_0{\bf P}'(w)={\bf
P}'(x_0 w)$ for $x=e,f$. Since $w$ is an arbitrary given word, we
conclude that $w$ is $\frak{gl}_{\mathcal{N}(m)}$-equivalent to
${\bf P}'(w)$. \qed

From the arguments in Corollary \ref{equivalence for nonstandard},
we can deduce the following, which is a generalization of
Proposition \ref{decomposition of Blambda}.
\begin{prop}\label{branching into kites}
For $\lambda\in\cP$ and $(\mu,\alpha)\in\C(m)$, the multiplicity of
$\B(\mu,\alpha)$ in $\B_{\N(m)}(\lambda)$ is equal to the number of
standard tableaux $T$ of shape $\lambda/\mu$ such that
$\alpha(T)=\alpha$ and the smallest entry $1$ of $T$ lies in the
first column of $\lambda$. In particular, $\B_{\N(m)}(\lambda)$ is
connected if the length of $\lambda$ is at most $m$.
\end{prop}

Let $\lambda,\mu,\nu$ be partitions with length no more than $m$.
Let $LR^{\lambda}_{\mu\nu}$ be the set of Littlewood-Richardson
tableaux of shape $\lambda/\mu$ with content $\nu$
(cf.~\cite{Fu,Mac}). We may also regard $LR^{\lambda}_{\mu\nu}$ as
the set of tableaux $S\in \B_{[m|0]}(\nu)$ such that $H_{\mu}\otimes
S$ is $\frak{gl}_{[m|0]}$-equivalent to $H_{\lambda}$ (see
\cite{Na}).

Let $(\lambda,\alpha),(\mu,\beta),(\nu,\gamma)\in\C(m)$ be given.
Define $\widetilde{LR}^{(\lambda,\alpha)}_{(\mu,\beta)(\nu,\gamma)}$
to be the set of quadruple $(S,T_1,T_2,w)$ satisfying the following
conditions;
\begin{itemize}
\item[(1)] $S\in LR^{\lambda}_{\eta\zeta}$, for some $\eta\subset
\mu$ and $\zeta\subset\nu$,

\item[(2)] $T_1\in
ST(\mu/\eta)$, and $T_2\in ST(\nu/\zeta)$,

\item[(3)] $w\in\mathcal{L}(w_{\alpha(T_1)\cdot\beta}\cup w_{\alpha(T_2)\cdot
\gamma})$, and the composition corresponding to the  descent set of
$w$ is $\alpha$ (see Section 4.1 and 4.2).

\item[(4)] Let $w_*$ be the smallest element in $w$, which is a node in $(\mu,\beta)$ or $(\nu,\gamma)$.
Then there is at least one $-1$ preceding the entry of $w_*$ when we
read the word associated to an element $T_1\otimes T_2$ in
$\B(\mu,\beta)\otimes \B(\nu,\gamma)$, where $\eta$ in $T_1$ and
$\zeta$ in $T_2$ are filled with $H_{\eta}$ and $S$, respectively.
\end{itemize}
Recall that given two compositions
$\sigma=(\sigma_1,\ldots,\sigma_r)$ and
$\tau=(\tau_1,\ldots,\tau_s)$, we mean by $\sigma\cdot\tau$ the
concatenation
$\sigma\cdot\tau=(\sigma_1,\ldots,\sigma_r,\tau_1,\ldots,\tau_s)$,
and given $T_1\otimes T_2\in \B(\mu,\beta)\otimes \B(\nu,\gamma)$,
the associated word is given by the juxtaposition $w(T_1)\cdot
w(T_2)$.

\begin{prop}\label{tensor product for Blambdaalpha}
Let $(\lambda,\alpha),(\mu,\beta),(\nu,\gamma)\in\C(m)$ be given.
The multiplicity of $\B(\lambda,\alpha)$ in $\B(\mu,\beta)\otimes
\B(\nu,\gamma)$ is equal to $|
\widetilde{LR}^{(\lambda,\alpha)}_{(\mu,\beta)(\nu,\gamma)}|$.
\end{prop}
\pf Let $U_1\otimes U_2\in \B(\mu,\alpha)\otimes \B(\nu,\beta)$ be
such that $e_{i}(U_1\otimes U_2)={\bf 0}$ for all $i\in I(m)$. Since
$U_1\otimes U_2$ is $\frak{gl}_{[m|0]}$-equivalent to $H_{\lambda}$
for some $\lambda\in \cP$, the subtableau of $U_1$ (resp. $U_2$)
consisting of negative entries is equal to $H_{\eta}$ (resp. $S\in
LR^{\lambda}_{\eta\zeta}$) for some $\eta\subset \mu$ and
$\zeta\subset \nu$. Next, the subtableau of $U_1$ (resp. $U_2$)
consisting of positive entries is $\frak{gl}_{\N}$-equivalent to a
unique quasi-ribbon tableau $U_1^+$ (resp. $U_2^+$) of shape
$\alpha(T_1)\cdot \beta$ (resp. $\alpha(T_2)\cdot \gamma$), where
$T_i$ ($i=1,2$) is the standard tableau uniquely determined by the
connected component of the subtableau in the body of $U_i$
consisting of positive entries (cf. Proposition \ref{decomposition
of Blambda}). Since $e_r (U_1\otimes U_2)={\bf 0}$ for all
$r\in\hf\mathbb{Z}_{> 0}$, $U_1^+\otimes U_2^+$ uniquely determines
a linear extension $w$ of $w_{\alpha(T_1)\cdot\beta}\cup
w_{\alpha(T_2)\cdot \gamma}$ by Proposition \ref{LR rule}. Then, we
can check that $w$ satisfies the condition (4) in the above
definition since $e_0 (U_1\otimes U_2)={\bf 0}$. Hence, we have
$(S,T_1,T_2,w)\in
\widetilde{LR}^{(\lambda,\alpha)}_{(\mu,\beta)(\nu,\gamma)}$, where
$\alpha$ is the composition of the descent set of $w$.

It is not difficult to see that the correspondence from $U_1\otimes
U_2$ to $(S,T_1,T_2,w)$ is reversible. Therefore, the number of
highest weight vectors in $\B(\mu,\beta)\otimes \B(\nu,\gamma)$,
whose connected components are isomorphic to $\B(\lambda,\alpha)$
for $(\lambda,\alpha)\in\C(m)$ is equal to $|
\widetilde{LR}^{(\lambda,\alpha)}_{(\mu,\beta)(\nu,\gamma)}|$. \qed

\subsection{Branching rule}
Suppose that ${\cS}$ and $\mathcal{T}$ are linearly ordered
$\mathbb{Z}_2$-graded sets. We say that ${\cS}$ and $\mathcal{T}$
are isomorphic if there exists a bijection from ${\cS}$ to
$\mathcal{T}$, which preserves both orderings and
$\mathbb{Z}_2$-gradings, and write ${\cS}\simeq \mathcal{T}$. Let
$\prec$ denote the linear ordering on $\cS$ and $\sigma$ a
permutation on ${\cS}$. Then $\sigma$ induces another linear
ordering $\prec_{\sigma}$ on ${\cS}$ given by $a\prec_{\sigma} b$ if
and only if $\sigma^{-1}(a)\prec\sigma^{-1}(b)$ for $a,b\in {\cS}$.
We denote by ${\cS}^{\sigma}$ the set ${\cS}$ with this new
ordering.

Suppose that ${\cS}=[m|n]$ for $m,n\geq 1$.  The associated Borel
subalgebra of $\frak{gl}_{[m|n]}$ is called {\it standard}. We
assume that $m\geq n$ for convenience. Let $\sigma$ be a permutation
on $[m|n]$. Note that $\sigma$ induces a natural isomorphism of Lie
superalgebras from $\frak{gl}_{[m|n]}$ to
$\frak{gl}_{[m|n]^{\sigma}}$, but the corresponding Borel
subalgebras are not conjugate in general. Moreover,
$\Delta_{[m|n]^{\sigma}}$ may have more than one odd isotropic
simple roots, while $\Delta_{[m|n]}$ has only one.

Now, let $\omega$ be a permutation on $[m|n]$ such that the number
of odd isotropic simple roots is maximal. We only have to consider a
shuffle of $\{-m,\ldots,-1\}$ and $\{1,\ldots,n\}$, which in fact
corresponds to a composite of a sequence of simple odd reflections
\cite{PS}. We may choose a unique $\omega$ such that
$[m|n]^{\omega}\simeq \N(p)^{\preccurlyeq q}$ with $p=m-n$ and
$q=n$.

\begin{ex}{\rm
\begin{equation*}
[4|2]^{\omega}=\{\,-4\prec_{\omega} -3\prec_{\omega} 1
\prec_{\omega} -2 \prec_{\omega} 2 \prec_{\omega} -1\,\}\simeq
\N(2)^{\preccurlyeq 2}.
\end{equation*}  }
\end{ex}


Then we have the following branching decomposition of
$\frak{gl}_{[m|n]}$-crystals into
$\frak{gl}_{[m|n]^{\omega}}$-crystals.

\begin{prop}\label{branching rule} For $\lambda\in\cP_{m|n}$,
there exists a $\frak{gl}_{[m|n]^{\omega}}$-crystal structure on
$\B_{[m|n]}(\lambda)$. For $(\mu,\alpha)\in \C(m-n)$, the
multiplicity of $\B_{[m|n]^{\omega}}(\mu,\alpha)$ in
$\B_{[m|n]}(\lambda)$ is equal to the number of standard tableaux
$T$ of shape $\lambda/\mu$ such that $\alpha(T)=\alpha$ and the
smallest entry $1$ of $T$ lies in the first column of $\lambda$.
\end{prop}
\pf Consider the bijection $\phi : \B_{[m|n]}(\lambda) \rightarrow
\B_{[m|n]^{\omega}}(\lambda)$  given by the super-analogue of
switching algorithm of Benkart, Sottile and Stroomer \cite{BKK} (see
also \cite{Hai89,K08,Rem90}). Given $T\in \B_{[m|n]}(\lambda)$, we
define $$x_{\alpha}T = \phi^{-1}(x_{\alpha} \phi(T))$$ for $x=e,f$
and $\alpha\in\Delta_{[m|n]^{\omega}}$. Hence, $\B_{[m|n]}(\lambda)$
becomes a $\frak{gl}_{[m|n]^{\omega}}$-crystal, which is isomorphic
to $\B_{[m|n]^{\omega}}(\lambda)$. We obtain the decomposition of
$\B_{[m|n]}(\lambda)$ by Proposition \ref{branching into kites}.
\qed

\begin{rem}{\rm
For a permutation $\sigma$ on $[m|n]$, $\B_{[m|n]}(\lambda)$ has a
$\frak{gl}_{[m|n]^{\sigma}}$-crystal structure by the same arguments
in Theorem \ref{branching rule}.
}
\end{rem}

\section{Super quasi-symmetric functions}

\subsection{Characters of $\frak{gl}_{\N(m)}$-crystals} Let ${\cS}$ be a
linearly ordered $\mathbb{Z}_2$-graded set. Let ${\bf z}={\bf
z}_{\cS}=\{\,z_b\,|\, b\in {\cS}\,\}$ be the set of formal
variables. For $\mu=\sum_{b\in {\cS}}\mu_b\epsilon_b\in P$, we set
${\bf z}^{\mu}=\prod_{b\in {\cS}}z_b^{\mu_b}$.  For a connected
component $C$ in $\W_{\cS}$, let ${\rm ch}C=\sum_{w\in C}{\bf
z}^{{\rm wt}(w)}$ be the character of $C$, which is a well-defined
formal series in ${\bf z}$. Define ${R}_{\cS}$ to be the
$\mathbb{Z}$-span of $\{\,{\rm ch}C(w)\,|\,w\in \W_{\cS}\,\}$.

\begin{prop}\label{chracter ring}
 $R_{\N(m)}$ is a commutative ring with a
$\mathbb{Z}$-basis $\{\,{\rm
ch}\B_{\N(m)}(\lambda,\alpha)\,|\,(\lambda,\alpha)\in\C(m)\,\}$ for
$m\in\mathbb{Z}_{\geq 0}$.
\end{prop}
\pf It follows from Corollary \ref{equivalence for nonstandard} and
Proposition \ref{tensor product for Blambdaalpha}. \qed

\begin{rem}{\rm

(1) The structure constants for $R_{\N(m)}$ are given in Proposition
\ref{tensor product for Blambdaalpha}.

(2) When $m=0$, $R_{\N}$ is isomorphic to the ring of
quasi-symmetric functions by Proposition \ref{LR rule} and
Proposition \ref{chracter ring}. In fact, ${\rm ch}\B_{\N}(\alpha)$
for $\alpha\in\C$ is a super quasi-symmetric function introduced in
\cite{HHRU}. }
\end{rem}

For $m\geq n\geq 1$, consider $R_{[m|n]}$ the ring of super
symmetric polynomials with $m+n$ variables which is spanned by hook
Schur polynomials
\begin{equation}
{\rm
ch}\B_{[m|n]}(\lambda)=hs_{\lambda}(z_{-m},\ldots,z_{-1};z_1,\ldots,z_n)
\end{equation}
for $\lambda\in\cP_{m|n}$ (cf.~\cite{Mac}). Note that ${\rm
ch}\B_{[m|n]}(\lambda)$ does not depend on the ordering of the
variables $z_{-m},\ldots,z_{-1},z_1,\ldots,z_n$, but the crystal
structure on $\B_{[m|n]}(\lambda)$ depends on the ordering on
$[m|n]$. Hence, we have a natural quasi-analogue for $R_{[m|n]}$ by
Proposition \ref{branching rule}.

\begin{prop} $R_{[m|n]}$ is a subring of
$R_{[m|n]^{\omega}}$, where $R_{[m|n]^{\omega}}$ is isomorphic to
$R_{\N(m-n)^{\preccurlyeq n}}$.
\end{prop}

Let $p\in\mathbb{Z}_{\geq 0}$ and $q\in\mathbb{Z}_{>0}$ be given. By
\eqref{elementsinzigzag}, we have
\begin{equation}
{\rm ch}\B_{\N^{\preccurlyeq q}}(2^{2q-1},1)=\prod_{i\in
\N^{\preccurlyeq q-\hf}}(z_{i}+z_{i+\hf}).
\end{equation}
In general, we have the following factorization property of ${\rm
ch}\B_{\N(p)^{\preccurlyeq q}}(\lambda,\alpha)$ (cf.~\cite{BR}).
\begin{prop}\label{factorization}
For $(\lambda,\alpha)\in \C(p)$, suppose that the number of corners
in $\alpha$ is either $2q-1$ or $2q$. Then
\begin{equation*}
\begin{split}
{\rm ch}\B_{\N(p)^{\preccurlyeq
q}}(\lambda,\alpha)&=hs_{\lambda}(z_{-p},\ldots,z_{-1};z_{\hf}){\rm
ch}\B_{\N^{\preccurlyeq q}}(\alpha) \\
&={\bf z}^{\mu}hs_{\lambda}(z_{-p},\ldots,z_{-1};z_{\hf})\prod_{i\in
\N^{\preccurlyeq q-\hf}}(z_{i}+z_{i+\hf}),
\end{split}
\end{equation*}
where $\mu={\rm wt}(H_{\alpha})-{\rm wt}(H_{(2^{2q-1},1)})$.
\end{prop}
\pf By Proposition \ref{stablity}, we observe that for $T\in
\B_{\N(p)^{\preccurlyeq q}}(\lambda,\alpha)$, the entries in $T$,
which are greater than $\hf$ and different from those at the same
place in $H_{(\lambda,\alpha)}$, occur only in the corners of its
tail. Hence we obtain a bijection from $\B_{\N(p)^{\preccurlyeq
q}}(\lambda,\alpha)$ to $\B_{\N(p)^{\preccurlyeq
\hf}}(\lambda)\times\B_{\N^{\preccurlyeq q}}(\alpha)$. Since ${\rm
ch}\B_{\N(p)^{\preccurlyeq \hf}}(\lambda)$ is a hook Schur
polynomial, this establishes the above identities. \qed

\subsection{Characterization of super quasi-symmetric functions}
We will give an algebraic characterization of $R_{\N(m)}$ or
$R_{\N(m)^{\preccurlyeq n}}$ for $m\in\mathbb{Z}_{\geq 0}$ and $n\in
\hf\mathbb{Z}_{>0}$, which is a quasi-analogue of Stembridge's
result on super symmetric polynomials \cite{St85}.

\begin{prop}[cf.~\cite{St85}]\label{ralpha}
Let $t$ be an indeterminate. For $(\lambda,\alpha)\in\C(m)$, we have
\begin{itemize}
\item[(1)] ${\rm ch}\B_{\N(m)}(\lambda,\alpha)$ is symmetric with respect to
$\{\,z_i\,|\,i=-m,\ldots,-1\,\}$,

\item[(2)] ${\rm ch}\B_{\N(m)}(\lambda,\alpha)|_{z_r=-z_{s}=t}$ is independent
of $t$, for $(r,s)=(-1,\hf)$ or $s=r+\hf$ with $r \in
\hf\mathbb{Z}_{>0}$. In particular, we have ${\rm
ch}\B_{\N(m)}(\lambda,\alpha)|_{z_r=-z_{r+\hf}=t}={\rm
ch}\B_{\N(m)\setminus\{r,r+\hf\}}(\lambda,\alpha)$ for $r \in
\hf\mathbb{Z}_{>0}$.
\end{itemize}
\end{prop}
\pf Note that when restricted to a $\frak{gl}_{\N(m)^{\preccurlyeq
\hf}}$-crystal, $\B_{\N(m)}(\lambda,\alpha)$ is a direct sum of
$\B_{\N(m)^{\preccurlyeq \hf}}(\mu)$'s for $\mu\in\cP_{m|1}$ up to
isomorphism. Hence the condition (1) is satisfied, and ${\rm
ch}\B_{\N(m)}(\lambda,\alpha)|_{z_{-1}=-z_{\hf}=t}$ is independent
of $t$ by the characterization of hook Schur polynomial in
\cite{St85} (see also \cite{Mac}).

Now, suppose that $r\in\hf\mathbb{Z}_{>0}$ is given. Then
$\B_{\N(m)}(\lambda,\alpha)$ is a disjoint union of
\begin{align*}
\B_{\N(m)}(\lambda,\alpha)^0&=\{\,T\,|\,e_rT=f_rT={\bf 0}\,\},\\
\B_{\N(m)}(\lambda,\alpha)^+&=\{\,T\,|\,f_rT\neq {\bf 0}\,\}, \\
\B_{\N(m)}(\lambda,\alpha)^-&=\{\,T\,|\,e_rT\neq {\bf 0}\,\}.
\end{align*}
Moreover, $f_r$ gives a bijection from
$\B_{\N(m)}(\lambda,\alpha)^+$ to $\B_{\N(m)}(\lambda,\alpha)^-$.
Since the number of occurrences of $r$ in $T\in
\B_{\N(m)}(\lambda,\alpha)^+$ has the different parity from that of
$f_rT\in \B_{\N(m)}(\lambda,\alpha)^-$. It follows that
\begin{equation*}
{\rm ch}\B_{\N(m)}(\lambda,\alpha)^+|_{z_r=-z_{r+\hf}=t}=-{\rm
ch}\B_{\N(m)}(\lambda,\alpha)^-|_{z_r=-z_{r+\hf}=t}.
\end{equation*}
This implies that ${\rm
ch}\B_{\N(m)}(\lambda,\alpha)|_{z_r=-z_{r+\hf}=t}={\rm
ch}\B_{\N(m)}(\lambda,\alpha)^0$ which is independent of $t$. On the
other hand, $\B_{\N(m)}(\lambda,\alpha)^0$ is the set of
semistandard tableaux of shape $(\lambda,\alpha)$ with entries in
$\N(m)\setminus \{r,r+\hf\}$. Hence, by definition, we have $${\rm
ch}\B_{\N(m)}(\lambda,\alpha)^0={\rm
ch}\B_{\N(m)\setminus\{r,r+\hf\}}(\lambda,\alpha).$$ The condition
(2) is satisfied. \qed

\begin{thm}
Suppose that $f$ is a polynomial in ${\bf z}_{\N(m)^{\preccurlyeq
n}}$ with integral coefficients. Then $f\in R_{\N(m)^{\preccurlyeq
n}}$ if and only if
\begin{itemize}
\item[(1)] $f$ is symmetric with respect to
$\{\,z_i\,|\,i=-m,\ldots,-1 \,\}$,

\item[(2)] $f|_{z_r=-z_{s}=t}$ is independent
of $t$, for $(r,s)=(-1,\hf)$ or $s=r+\hf$ with $r \in
\hf\mathbb{Z}_{>0}$.
\end{itemize}
\end{thm}
\pf Let $\mathscr{R}_{m,n}$ be the ring of polynomials $f$ in ${\bf
z}_{\N(m)^{\preccurlyeq n}}$ with integral coefficients satisfying
(1) and (2). By Proposition \ref{ralpha}, it is enough to show that
$\mathscr{R}_{m,n}\subset R_{\N(m)^{\preccurlyeq n}}$. We will use
induction on $n$. It is clear when $n=\hf$ by the characterization
of super symmetric polynomials \cite{St85}.

Suppose that $n\succcurlyeq 1$ and $f\in \mathscr{R}_{m,n}$ is
given. By induction hypothesis, we have
\begin{equation*}
f|_{z_{n-\hf}=-z_{n}=t}=\sum_{(\lambda,\alpha)\in\C(m)}c_{(\lambda,\alpha)}
g_{(\lambda,\alpha)},
\end{equation*}
where $g_{(\lambda,\alpha)}={\rm ch}\B_{\N(m)^{\preccurlyeq
n-\hf}}(\lambda,\alpha)$ with $c_{(\lambda,\alpha)}\in\mathbb{Z}$.
Put
\begin{equation*}
h=f-\sum_{(\lambda,\alpha)}c_{(\lambda,\alpha)}{\rm
ch}\B_{\N(m)^{\preccurlyeq n}}(\lambda,\alpha).
\end{equation*}
Since $h|_{z_{n-\hf}=-z_{n}=t}=0$, we have
$h=(z_{n-\hf}+z_{n})h^{(1)}$ for some polynomial $h^{(1)}$. Consider
\begin{equation*}
h|_{z_{n-1}=-z_{n-\hf}=t}=(-t+z_{n-1})\left(h^{(1)}|_{z_{n-1}=-z_{n-\hf}=t}\right).
\end{equation*}
Since the left-hand side is independent of $t$, the right-hand side
must be zero, and hence $(z_{n-1}+z_{n-\hf})$ divides $h^{(1)}$.
Repeating this procedure, we conclude that $\prod_{i\in
\N^{\preccurlyeq n-\hf}}(z_{i}+z_{i+\hf})$ divides $h$. Moreover,
since $h$ is super symmetric, we have $h=\prod_{i\in
\N^{\preccurlyeq n-\hf}}(z_{i}+z_{i+\hf})
\prod_{i=-m}^{-1}(z_{i}+z_{\hf})\ k$, where $k$ is a symmetric
polynomial in $\{\,z_{-m},\ldots,z_{-1}\,\}$ with coefficients in
$\mathbb{Z}[z_{\hf},\ldots,z_{n}]$. By Proposition
\ref{factorization}, it follows that $h$ is an integral linear
combination of ${\rm ch}\B_{\N(m)^{\preccurlyeq
n}}(\lambda,\alpha)$'s and so is $f$. This completes the induction.
\qed

\begin{rem}{\rm
(1) For each $k$, let $R^k_{\N(m)^{\preccurlyeq n}}$ be the subgroup
of polynomials of homogeneous degree $k$ in $R_{\N(m)^{\preccurlyeq
n}}$, whose inductive limit is defined in a standard way and denoted
by $R^k_{\N(m)}$. Then we have $R_{\N(m)}=\bigoplus_{k\geq
0}R^k_{\N(m)}$.

(2) Let $R_{\N}^{\circ}$ be the subring of $R_{\N}$ with a
$\mathbb{Z}$-basis $\{\,{\rm
ch}\B_{\N}(\lambda)\,|\,\lambda\in\cP\,\}$. It is well-known as the
ring of super symmetric functions \cite{BR,Mac}. Then using the
characterization of super symmetric polynomials in \cite{St85}, one
can easily deduce that for $f\in R_{\N}$, $f$ is super symmetric,
that is, $f\in R_{\N}^{\circ}$ if and only if $f$ is symmetric with
respect to both ${\bf z}_{\mathbb{N}}$ and ${\bf
z}_{\hf+\mathbb{Z}_{\geq 0}}$.}
\end{rem}

{\bf Acknowledgement} Part of this this work was done during the
author's visit at National Taiwan University 2006 Summer. He thanks
Prof. S.-J. Cheng for the invitation and many helpful discussions.
He also thanks Prof. S.-J. Kang and G. Benkart for their interests
in this work.

\end{document}